\documentclass[10pt]{article}

\usepackage{amssymb}
\usepackage{amsmath}
\usepackage{graphicx}								% to make \includegraphics work
\usepackage{subfig}
\usepackage{amsfonts} 							% to make \mathbb work
\usepackage{dsfont}
\usepackage{bigints} 								% to use big integrals 
\usepackage{algorithm}              % algorithm environment
\usepackage{algorithmic}
\usepackage{booktabs} 							% correct spacing in tables
\usepackage{multirow}								% multirows in tables
\usepackage{setspace}
\usepackage{mathrsfs}

\usepackage{geometry}
\usepackage{tikz}	
\usetikzlibrary{arrows,chains,matrix,positioning,scopes}

\makeatletter
\tikzset{join/.code=\tikzset{after node path={%
\ifx\tikzchainprevious\pgfutil@empty\else(\tikzchainprevious)% 
edge[every join]#1(\tikzchaincurrent)\fi}}}

\makeatother
\tikzset{>=stealth',every on chain/.append style={join}, every join/.style={->}}
\tikzstyle{labeled}=[execute at begin node=$\scriptstyle,   execute at end node=$]

\newcommand {\Int}   {\int\limits}
\newcommand {\Sum}   {\sum\limits}

\newcommand {\IntQT} {\Int_{Q_T}}

\newcommand {\IntT}  {\Int_{t_k}^{t_{k+1}}}
\newcommand {\IntO}  {\Int_\Omega}
\newcommand {\IntGammaN}  {\Int_{\Gamma_N}}

\newcommand {\CF}    {C_{\mathrm{F\Omega}}}
\newcommand {\Ctr}   {C_{\mathrm{tr}}}
\newcommand {\Rd}    {{\mathds{R}}^d}

\newcommand {\R}     {\mathscr{R}}
\newcommand {\I}     {\mathscr{I}}
\newcommand {\Ieff}  {{I_{\rm eff}}}
\newcommand {\Ieffmaj}  {{I_{\rm eff}^{\overline{\mathrm M}}}}
\newcommand {\Ieffmin}  {{I_{\rm eff}^{\underline{\mathrm M}}}}
\newcommand {\Ieffmajzero}  {{I_{\rm eff}^{\overline{\mathrm M}_{(0)}}}}
\newcommand {\Ieffmajmuhat}  {{I_{\rm eff}^{\overline{\mathrm M}_{(\hat \mu)}}}}
\newcommand {\Ieffmajone}  {{I_{\rm eff}^{\overline{\mathrm M}_{(1)}}}}

\def \NormA#1  {{\mid\!\mid\!\mid\! #1 \!\mid\!\mid\!\mid}^2 }   % |||...|||
\def \NormAinverse#1  { {\mid\!\mid\!\mid\! #1 \!\mid\!\mid\!\mid}^2_* }   % |||...|||_*
\def \NormQT#1 {{\mid\!\mid\!\mid #1 \mid\!\mid\!\mid}^2_{Q_T}}   % ||...||
\def \NormQO#1 {{\mid\!\mid\!\mid #1 \mid\!\mid\!\mid}^2_{Q^0}}   % ||...||
\def \NormQk#1 {{\mid\!\mid\!\mid #1 \mid\!\mid\!\mid}^2_{Q^k}}   % ||...||
\def \Normf#1  {\Big \lceil #1 \Big\rceil_\Omega}
\def \dvrg     {\mathrm{div} \:}	

\def \dt       {\mathrm{\:d}t}
\def \dx       {\mathrm{\:d}x}
\def \dxt      {\mathrm{\:d}x\mathrm{d}t}
\def \dst      {\mathrm{\:d}s\mathrm{d}t}
\def \ds       {\mathrm{\:d}s}
\def \Normt#1  {\mid\!\mid\!\mid #1 \mid\!\mid\!\mid}   % |||...|||

\def\L#1{L^{#1}}
\def\H#1{H^{#1}}
\def\Ho#1{{\mathring{H}}^{#1}}

\def\Maj{\overline{\mathrm M}^2_{(\delta, \, \gamma, \, \mu)}}
\def\Majmu{\overline{\mathrm M}^2_{(\mu )}}
\def\Majmuhat{\overline{\mathrm M}^2_{(\hat \mu )}}
\def\Majone{\overline{\mathrm M}^2_{(1)}}
\def\Majzero{\overline{\mathrm M}^2_{(0)}}

\def\maj{\overline{\mathrm M}^{2}}
\def\majk#1{\overline{\mathrm M}^{2\,#1}}

\def\Min{{\underline{\mathrm M}^2}}
\def\mink#1{\underline{\mathrm M}^{2\,#1}}
\def\error{{[e]\,}^2}
\def\mdI{\overline{\mathrm m}^2_{\mathrm{d}}}
\def\mfI{\overline{\mathrm m}^2_{\mathrm{f}}}
\def\Marker{\mathbb M}

\newcommand{\minus}{\scalebox{0.5}[1.0]{\( - \)}}

\newtheorem{thm}{Theorem}

\newtheorem{example}{Example}
\newtheorem{rmk}{Remark}

\title{{\bf Computable estimates of the distance to the exact solution of the 
evolutionary reaction-diffusion equation}}
\author{
Svetlana Matculevich\\
\small {\rm svetlana.v.matculevich@jyu.fi}, \\
\small {\it Dept. of Mathematical Information Technology, 
C321.4, Agora, P.O. Box 35, 
FI-40014 University of Jyv{\"a}skyl{\"a}, Finland}
\\\\
Sergey Repin\\
\small {\rm repin@pdmi.ras.ru}, \\
\small {\it V.A. Steklov Institute of Mathematics at St. Petersburg, 
191011, Fontanka 27, St.Petersburg, Russia}
}

\begin{document}

\maketitle
%\begin{frontmatter}


%% use the cortextf command for the associated footnote;

%\author[focal]{S.~Repin}
%\ead{repin@pdmi.ras.ru}
%\address[focal]{V.A. Steklov Institute of Mathematics at St. Petersburg, 191011, Fontanka 27, 
%St.Petersburg, Russia}

%--------------------------------------------------------------------------------------%

\begin{abstract}

We derive guaranteed bounds of distance to the exact solution of the evolutionary 
reaction-diffusion problem with mixed Dirichlet--Neumann boundary condition. It is shown 
that two-sided error estimates are directly computable and equivalent to the error. 
Numerical experiments confirm that estimates provide accurate two-sided bounds of the 
overall error and generate efficient indicators of local error distribution.

\end{abstract}

%------------------------------------------------------------------------------------ -%
\section{Problem statement}
%--------------------------------------------------------------------------------------%

%--------------------------------------------------------------------------------------%
Let $\Omega \in \Rd$ ($d = $ 1, 2, or 3)  be a bounded connected domain with Lipchitz 
continuous boundary $\partial \Omega$, which consists of two measurable 
non-intersecting parts $\Gamma_D$ and $\Gamma_N$ associated with Dirichlet and Neumann 
boundary conditions, respectively. Let $Q_T$ denote the space-time cylinder 
$Q_T := \Omega \times (0, T)$, $T > 0$. The lateral surface of $Q_T$ 
is denoted by $S_T := \partial\Omega \times [0, T] = 
\big( \Gamma_D \cup \Gamma_N \big) \times [0, T] = 
S_D \cup S_N$. 

%--------------------------------------------------------------------------------------%

We consider the classical reaction-diffusion initial boundary value problem
\begin{alignat}{3}
   u_t - \nabla \cdot p + \lambda u & =\, f,	      & \quad (x, t) \in Q_T,\label{eq:parabolic-problem-equation}\\
                                  p & =\, A \nabla u, & \quad (x, t) \in Q_T,\nonumber\\
                      			u(x, 0) & =\, \varphi,		& \quad x \in \Omega,\label{eq:parabolic-problem-initial-condition}\\
                      			      u & =\, 0,		      & \quad (x, t) \in S_D,\label{eq:parabolic-problem-dirichlet-boundary-condition}\\
   			                  p \cdot n & =\, g,					& \quad (x, t) \in S_N,\label{eq:parabolic-problem-robin-boundary-condition}
\end{alignat}
where $n$ denotes the vector of unit outward normal to $\partial\Omega$, and
\begin{equation}
f(x, t) \in \L{2}(Q_T), \quad  
\varphi(x) \in \L{2}(\Omega), \quad 
{g(x, t) \in \L{2}\left(0, T; \L{2} (\Gamma_N)\right).}
\label{eq:problem-condition}
\end{equation}
The function $\lambda$ entering the reaction part of 
(\ref{eq:parabolic-problem-equation}) is non-negative bounded function, 
and its values may 
drastically vary on different parts of the domain. Also, we assume that for a.e. 
$t \in (0, T)$ the matrix $A$ is symmetric and satisfies the condition
\begin{equation}
\nu_1 |\xi|^2 \leq A(x) \: \xi \cdot \xi \leq \nu_2 |\xi|^2,\quad
\xi \in \Rd,\quad 0 < \nu_1 \leq \nu_2 < \infty, \quad
\mbox{for a.e.} \quad x \in \Omega.
\label{eq:operator-a}
\end{equation}
Henceforth, we use the following notation
\begin{equation}
  \| \: \tau \: \|^2_A = \Int_\Omega A \tau \cdot \tau \dx, \quad 	
  \| \: \tau \: \|^2_{A^{-1}} = \Int_\Omega A^{-1} \tau \cdot \tau \dx, \quad
	\NormA{\tau} =  \Int_{Q_T} A \tau \cdot \tau \dxt, \quad
	\NormAinverse{\tau} =  \Int_{Q_T} A^{-1} \tau \cdot \tau \dxt.
	\label{eq:norms}
\end{equation}
%
%--------------------------------------------------------------------------------------%
%
By $\| \cdot \|_\Omega$ and $\| \cdot \|_{Q_T}$ we denote the norms in  $\L{2}(\Omega)$ and 
$\L{2}(Q_T)$, respectively. The space of functions $g(x, t)$ with the norm 
$\Int_0^T \| g(\cdot, t)\|_{\Omega} \dt$ is denoted by $\L{2, 1} (Q_T)$. 
%
%--------------------------------------------------------------------------------------%
%
$\Ho{1}(Q_T)$ is a subspace of $\H{1}(Q_T)$, which contains the functions satisfying 
(\ref{eq:parabolic-problem-dirichlet-boundary-condition}), 
$\Ho{1, 0} (Q_T) := \L{2}\Big( (0, T); \Ho{1}(\Omega) \Big)$, and
$V_2(Q_T) := \Ho{1, 0} (Q_T) \cap \L{\infty}\Big( (0, T); \L{2}(\Omega) \Big)$.
%
%--------------------------------------------------------------------------------------%
%
The space $V^{1, 0}_2 (Q_T) := \Ho{1, 0} (Q_T) \cap C \Big([0, T]; \L{2}(\Omega)\Big)$ 
is a subspace of $V_2(Q_T)$ and contains the functions with traces from 
$\L{2}(\Omega)$ for all $t \in [0, T]$, which continuously changing with 
$t \in [0, T]$ in the $\L{2}(\Omega)$ norm.  
%--------------------------------------------------------------------------------------%

The generalized statement of 
(\ref{eq:parabolic-problem-equation})--(\ref{eq:parabolic-problem-robin-boundary-condition}) 
is as follows: find a function $u(x, t) \in V^{1, 0}_2 (Q_T)$ satisfying the integral 
identity
\begin{multline}
	\Int_{\Omega} \Big( u(x, T) \eta(x, T) - u(x, 0) \eta(x, 0) \Big) \dx - 
	\Int_{Q_T} u \eta_t \dxt + 
	\Int_{Q_T} A \nabla{u} \cdot \nabla{\eta} \dxt +  
	\Int_{Q_T} \lambda u \eta \dxt = \\
	\Int_{Q_T} f \eta \dxt + \Int_{S_N} g \eta \dst, \quad \forall \eta \in \Ho{1}(Q_T).
	\label{eq:generalized-statement}
\end{multline}
%
%--------------------------------------------------------------------------------------%
Well-known classical solvability results (see, e.g., 
\cite{Ladyzhenskaya1985, Ladyzhenskayaetall1967, Evans2010}) guarantee that 
(\ref{eq:parabolic-problem-equation})--(\ref{eq:parabolic-problem-robin-boundary-condition}) 
has a unique solution in $V^{1, 0}_2(Q_T)$ provided that the conditions 
(\ref{eq:problem-condition}) hold.

%--------------------------------------------------------------------------------------%
Assume that $v \in \Ho{1}(Q_T)$ presents a certain
numerical approximation of $u$. Our goal is to deduce accurate and explicitly 
computable estimates of the distance between $u$ and $v$. For this purpose, we use the 
norm
\begin{equation}	
	[ u - v ]^2_{(\nu,\, \theta,\, \zeta)} := 
	\NormA{\nu\, \nabla (u - v) \,} \: + \:
	\| \, \theta \, (u - v) \, \|^2_{Q_T} + \:
	\| \, \zeta\, (u - v)(x, T) \, \|^2_{\Omega},
	\label{eq:energy-norm-for-reaction-diff-evolutionary-problem}
\end{equation}
where $\nu, \theta$, and $\zeta$ are certain positive weights (weight-functions), 
which balance three 
components of the error. They can be selected in different ways so that 
(\ref{eq:energy-norm-for-reaction-diff-evolutionary-problem}) presents a common form 
of a wide spectrum of different measures.
%--------------------------------------------------------------------------------------%

A fully computable and guaranteed upper bound of $[ u - v ]^2_{(\nu,\, \theta,\, \zeta)}$ 
is derived in Theorem \ref{th:theorem-minimum-of-majorant-I} with the help of the method 
originally introduced in \cite{Repin2002}. In \cite{RepinTomar2010}, the method was 
applied to problems with convection, and in \cite{NeittaanmakiRepin2010} the guaranteed error 
majorants were derived for the Stokes problem. In Section \ref{sec:majorant-I}, we 
combine this method with the technique suggested in \cite{RepinSauter2006} for the 
stationary reaction-diffusion problem, which makes it possible to obtain the efficient error 
majorants for problems with drastically different values of the reaction function. 
Theorem \ref{th:theorem-minimum-of-majorant-I} presents such an estimate for the problem 
(\ref{eq:parabolic-problem-equation})--(\ref{eq:parabolic-problem-robin-boundary-condition}). 
In Section \ref{sec:minorant}, we derive a guaranteed and fully computable lower bound 
of the error (\ref{eq:energy-norm-for-reaction-diff-evolutionary-problem}) (Theorem 
\ref{th:theorem-minimum-of-minorant}). Sections \ref{sec:incremental} and 
\ref{sec:numerics} are devoted to practical applications of the estimates. 
In them, we discuss numerical results obtained for several typical examples, which confirm the efficiency of two-sided bounds.

%--------------------------------------------------------------------------------------%
\section{Error majorant}
\label{sec:majorant-I}
%--------------------------------------------------------------------------------------%

%--------------------------------------------------------------------------------------%
Let $e(x, t) := (u - v)(x, t)$ denote the  deviation of $v\in \Ho{1}(Q_T)$ from the 
exact solution $u$. From (\ref{eq:generalized-statement}), it follows that
\begin{multline}
	\Int_{\Omega} \big( e(x, T) \eta(x, T) - e(x, 0) \eta(x, 0) \big) \dx -
	\Int_{Q_T} e \eta_t \dxt + 
	\Int_{Q_T} A \nabla{e} \cdot \nabla{\eta} \dxt + 
	\Int_{Q_T} \lambda e \eta \dxt = \\
	\Int_{Q_T} \left(f - v_t - \lambda v \right)\eta \dxt -
	\Int_{Q_T} A \nabla{v} \cdot \nabla{\eta} \dxt + \Int_{S_N} g \eta \dst.
	\label{eq:gen-stat-waveu-w}
\end{multline}
Since $e \in \Ho{1}(Q_T)$, we can set $\eta = e$ and, using the relation
\begin{equation}
	\Int_{\Omega} \left( e^2(x, T) - e^2(x, 0)\right) \dx - \Int_{Q_T} e e_t \dxt = 	
	\frac12 \left(\| \, e(x, T) \, \|^2_{\Omega} - \| \, e(x, 0) \,\|^2_{\Omega} \right),
\end{equation}
obtain
\begin{multline}
	\frac12 \| \, e(x, T) \,\|^2_{\Omega} + \NormA{\nabla{e}} + 
	\left \| \sqrt{\lambda}e \, \right \|^2_{Q_T} = \\
	\Int_{Q_T} \left(f - v_t - \lambda v \right) e \dxt -
	\Int_{Q_T} A \nabla{v} \cdot \nabla e \dxt + 
	\Int_{S_N} g e \dst + \frac12 \| \, e(x, 0) \, \|^2_{\Omega}\,.
	\label{eq:energy-balance-equation}
\end{multline}
This relation is a form of the `energy-balance' identity in terms of deviation. It 
plays an important role in subsequent analysis. Now, we introduce an additional 
variable $y \in Y_{\dvrg}(Q_T)$, where
\begin{equation}
	Y_{\dvrg}(Q_T) := \left \{ \: y(x,t) \in \L{2}(\Omega, \Rd) \;\; \big | \;\; 
	\dvrg y(x,t) \in \L{2}(\Omega), \;\;
	y \cdot n \in \L{2}(\Gamma_N) \;\; 
	{\mbox{for a.e.}} \;\; t \in (0, T) \: \right \}.
	\label{eq:y-set-div}
\end{equation}
%
%--------------------------------------------------------------------------------------%
%
\begin{thm}
\label{th:theorem-minimum-of-majorant-I}
(i) For any 
$v \in \Ho{1}(Q_T)$ and $y \in Y_{\dvrg}(Q_T)$ the following inequality holds
\begin{multline}
	(2 - \delta)\NormA{\nabla e} + 
	\left( 2 - \frac1\gamma \right) \left \| \sqrt{\lambda}e \, \right \|^2_{Q_T}  + 
	\| \, e (x, T) \, \|^2_{\Omega} =: 
	\error_{(\overline{\nu},\, \overline{\theta},\, \overline{\zeta})} \leq 
	\Maj (v, y) := \\
	\| \, e (0, x) \, \|^2_{\Omega} + 
	        \Int_0^T \Bigg ( \gamma \bigg\| \, \frac{\mu}{\sqrt{\lambda}} \, \R_f(v, y) \, \bigg\|^2_{\Omega} \: + 
	                 \alpha_1(t) \frac{\CF^2}{\nu_1} \| \, (1 - \mu) \, \R_f(v, y) \, \|^2_{\Omega} \: + \\
									 \alpha_2(t) \| \, \R_d(v, y) \,\|^2_{A^{-1}} +
									 \alpha_3(t) \frac{\Ctr^2}{\nu_1} \big\| \, \R_b (v, y) \, \big\|^2_{\Gamma_N} \Bigg ) \dt,
	\label{eq:majorant-1}
\end{multline}
where 
\begin{equation}
\R_f (v, y) = f - v_t - \lambda v + \dvrg y, \quad  
\R_d (v, y) = y - A \nabla{v}, \quad 
\R_b (v, y) = g - y \cdot n, 
\end{equation}
$\CF$ is the constant from the Friedrichs' inequality
\begin{equation}
\| \, \bar{\eta} \, \|_{\Omega} \leq \CF \| \, \nabla \bar{\eta} \, \|_{\Omega}, \qquad \forall \bar{\eta} \in \Ho{1}(\Omega), 
\end{equation}
$\Ctr$ is the constant in the trace inequality 
\begin{equation}
\| \, \bar{\eta} \, \|_{\Gamma_N} \leq \Ctr \| \, \nabla \bar{\eta} \, \|_{\Omega}, \qquad \forall \bar{\eta} \in \Ho{1}(\Omega).
\end{equation}
Here, 
$\overline{\nu} = \sqrt{2 - \delta}$, 
$\overline{\theta} = \sqrt{2 - \frac1\gamma}$, 
$\overline{\zeta} = 1$ are positive weights, where 
$\delta \in (0, 2]$, 
$\gamma \in \Big[\frac{1}{2}, + \infty\Big[$;
$\mu(x, t)$ is a real-valued function taking values in $[0, 1]$; and 
$\alpha_1(t)$, $\alpha_2(t)$, $\alpha_3(t)$ are positive scalar-valued functions 
satisfying the relation
\begin{equation}
	\frac{1}{\alpha_1(t)} + \frac{1}{\alpha_2(t)} +
	\frac{1}{\alpha_3(t)} = \delta.
	\label{eq:alpha-relation}
\end{equation}
%--------------------------------------------------------------------------------------%

\noindent
(ii) For any $\delta \in (0, 2]$, 
$\gamma \in \Big[\frac{1}{2}, + \infty\Big[$, and 
a real-valued function $\mu(x, t)$ taking values in $[0, 1]$, the lower bound
of the variation problem 
\begin{equation}
\inf\limits_{
\begin{array}{c}
v \in \Ho{1}(Q_T)\\
y \in Y_{\dvrg}(Q_T)
\end{array}
} \Maj (v, y)
\label{eq:inf-maj-I}
\end{equation}
is zero, and it is attained if and only if $v = u$ and $y = A \nabla u$.
\end{thm}

%--------------------------------------------------------------------------------------%
%
\noindent {\bf Proof.}
(i) We  transform the right-hand side of (\ref{eq:energy-balance-equation}) by 
means of the relation
\begin{equation}
	\Int_{Q_T} \dvrg y \, e \dxt + \Int_{Q_T} y \cdot \nabla{e}\dxt =
	\Int_{S_N} y \cdot n \,e \dst,
	\label{eq:y-rel}
\end{equation}
which yields
\begin{equation}
	\frac12 \| \, e (x, T) \, \|^2_{\Omega} + \NormA{\nabla{e}} + 
  \left \| \sqrt{\lambda}e \, \right \|^2_{Q_T} =
	\I_f + \I_d + \I_b + \frac12 \| \, e(x, 0) \, \|^2_{\Omega},
	\label{eq:gen-stat-u-v-norm-y}
\end{equation}
where
\begin{equation}
	\I_f = \Int_{Q_T} \R_f (v, y) \, e \dxt, \quad
	\I_d = \Int_{Q_T} \R_d (v, y) \cdot \nabla{e} \dxt, \quad 
	\I_b = \Int_{S_N} \R_b (v, y) \, e \dst. 
	\label{eq:Ir-Id-Ib-terms}
\end{equation}
%
%--------------------------------------------------------------------------------------%
By means of the H\"older's inequality, we find that 
\begin{equation}
	\I_d = \Int_{Q_T} \R_d (v, y) \cdot \nabla{e} \dxt \leq
	       \Int_0^T   \left \| \, \R_d (v, y) \, \right\|_{A^{-1}} \| \, \nabla{e} \, \|_A \dt
	\label{eq:id-estimate}
\end{equation}	
and
\begin{equation}
	\I_b = \Int_{S_N} \R_b (v, y) \, e \dst \leq 
	      \Int_0^T \left\| \, \R_b (v, y) \, \right\|_{\Gamma_N} \|\, e\, \|_{\Gamma_N} \dt \leq 
 	      \Int_0^T \left\| \, \R_b (v, y) \, \right\|_{\Gamma_N} \frac{\Ctr}{\sqrt{\nu_1}} \|\, \nabla{e}\, \|_A \dt,
	\label{eq:ib-estimate}
\end{equation}	
{where $\nu_1$ appears due to (\ref{eq:operator-a}). 
%in the right-hand side of (\ref{eq:ib-estimate}) due to the inequality 
%
%\begin{equation}
%	\| \, \xi \,\| \leq \frac{1}{\sqrt{\nu_1}} \| \,\xi \,\|_{A}, \quad 
%	\forall \xi \in \Ho{1} (\Omega).
%\end{equation}
}
Let $\mu(x, t)$ be a real-valued function taking values in $[0, 1]$. Then, we estimate 
the term $\I_f$ as follows:
\begin{equation}
	\I_f \leq \bigintsss\limits_0^T \left (
	         \bigg\| \,  \frac{\mu}{\sqrt{\lambda}} \R_f(v, y) \,  \bigg\|_{\Omega}
	         \left \|  \sqrt{\lambda}e  \,  \right \|_{\Omega} +
	         \frac{\CF}{\sqrt{\nu_1}} \left \| \,  (1 - \mu)\, \R_f(v, y) \,  \right \|_{\Omega}
	         \|\, \nabla e \, \|_A \right )\dt.
  \label{eq:if-estimate}
\end{equation}	
In \cite{RepinSauter2006}, this decomposition was used in order to overcome difficulties 
arising in the stationary problem  if $\lambda$ is small (or close to zero) in some 
parts of the domain (a more detailed study of this form of 
{the majorant} can be found in 
\cite{NeittaanmakiRepinMaxwell2010, MaliNeittaanmakiRepin2013}).

%--------------------------------------------------------------------------------------%

Combining (\ref{eq:id-estimate}), (\ref{eq:id-estimate}), and (\ref{eq:if-estimate}), 
we obtain
\begin{multline}
	\frac12 \| \, e(x, T) \, \|^2_{\Omega} + \NormA{\nabla e} + 
  \left \| \, \sqrt{\lambda}e \, \right \|^2_{Q_T} \leq 
	\frac12 \| \, e(x, 0) \, \|^2_{\Omega} + 
	\Int_0^T \Bigg (
	\bigg\| \, \frac{\mu}{\sqrt{\lambda}} \, \R_f(v, y) \, \bigg\|_{\Omega} \Big \| \, \sqrt{\lambda} \, e \Big \|_{\Omega} + \\
  \frac{\CF}{\sqrt{\nu_1}} \left \| \, (1 - \mu) \, \R_f(v, y) \, \right \|_{\Omega} \| \, \nabla e \, \|_A + 
	\left \| \, \R_d (v, y) \, \right \|_{A^{-1}} \| \, \nabla e \,\|_A + \left\| \, \R_b (v, y) \, \right\|_{\Gamma_N} \frac{\Ctr}{\sqrt{\nu_1}} \| \, \nabla e \, \|_A \Bigg ) \dt.
	\label{eq:estimate}
\end{multline} 	      	      
The second term on the right-hand side of (\ref{eq:estimate}) can be estimated by the
Young--Fenchel inequality
%--------------------------------------------------------------------------------------%
%
\begin{align}	
	\Int_0^T \bigg\| \,\frac{\mu}{\sqrt{\lambda}} \, \R_f(v, y) \,\bigg\|_{\Omega} 
		 \left \| \sqrt{\lambda} \, e \, \right \|_{\Omega} \dt \leq   
	\Int_0^T \left ( \frac\gamma2      \bigg\| \, \frac{\mu}{\sqrt{\lambda}} \, \R_f(v, y) \, \bigg\|^2_{\Omega} +
		\frac{1}{2\gamma} \left \| \, \sqrt{\lambda} e \, \right  \|^2_{\Omega} \right) \dt,
	\label{eq:young-fenchel-1}
\end{align}	
where $\gamma$ is an arbitrary real parameter from $\Big[\frac{1}{2}, + \infty\Big[$.
 Analogously,
\begin{align}	
	\Int_0^T \frac{\CF}{\sqrt{\nu_1}}\left \| \, (1 - \mu) \, \R_f (v, y)  \, \right \|_{\Omega} \|\, \nabla e \, \|_A \dt & \leq  	         
	\Int_0^T \left(\frac{\alpha_1(t)}{2}  \frac{\CF^2}{\nu_1} \left \| \, (1 - \mu) \, \R_f (v, y) \, \right \|^2_{\Omega} +
	               \frac{1}{2\alpha_1(t)} \| \, \nabla e \, \|^2_A \right) \dt, 
\end{align}	
	
\begin{align}
	\Int_0^T \| \, \R_d (v, y) \, \|_{A^{-1}} \|\, \nabla e \, \|_A  \dt & \leq 
  \Int_0^T \left( \frac{\alpha_2(t)}{2}  \left \| \, \R_d (v, y) \, \right \|^2_{A^{-1}} +
							    \frac{1}{2\alpha_2(t)} \|\, \nabla e \,\|^2_A \right) \dt,
\end{align}	
and
\begin{align}	
	\Int_0^T \left\| \, \R_b (v, y) \, \right\|_{\Gamma_N}
	               \frac{ \Ctr }{\sqrt{\nu_1}}\| \, \nabla e \, \|_A \dt \leq 
	\Int_0^T \left(\frac{\alpha_3(t)}{2} \frac{\Ctr^2}{\nu_1}  
					 \left \| \, \R_b (v, y)\, \right \|^2_{\Gamma_N} +
					 \frac{1}{2\alpha_3(t)} \| \, \nabla e \, \|^2_A \right) \dt.
	\label{eq:young-fenchel-4}
\end{align}
Here, $\alpha_1(t)$, $\alpha_2(t)$, and $\alpha_3(t)$ are functions satisfying 
(\ref{eq:alpha-relation}). The estimate (\ref{eq:majorant-1}) follows from
(\ref{eq:young-fenchel-1})--(\ref{eq:young-fenchel-4}).

%--------------------------------------------------------------------------------------%
(ii) Existence of the pair $(v, y) \in \Ho{1}(Q_T) \times Y_{\dvrg}(Q_T)$ minimizing 
the functional $\Maj (v, y)$ can be proven straightforwardly. Indeed, let $v = u$ and
$y = A \nabla u$. Since $\dvrg (A \nabla u) \in \L{2}(Q_T)$, we see that 
$y \in Y_{\dvrg}(Q_T)$. In this case, 
$e(0, x) = (u - v)(0, x) = \varphi(x) - v(0, x) = 0$, \,
$\R_f(u, A \nabla u) = f - u_t - \lambda u + \dvrg A \nabla u = 0$, \,
$\R_d (u, A \nabla u) = A \nabla u - A \nabla{u} = 0 $,  and 
$\R_b (v, y) = g - A \nabla u \cdot n = 0$.  Then, 
$\Maj ( u, A \nabla u ) = 0$, and 
the exact lower bound is attained. 

On the other hand, if $\Maj (v, y) = 0$, then $v$ satisfies the 
initial and boundary conditions, and for 
{ a.e.} $(x, t) \in Q_T$ the following 
relations hold:
\begin{equation}
  y = A \nabla v, \quad f - v_t - \lambda v + \dvrg y = 0.
	\label{eq:proof-relations-set-2}  
\end{equation}
From (\ref{eq:proof-relations-set-2}), it follows that
\begin{equation}
	\Int_{Q_T} ( f - v_t - \lambda v ) \eta \dxt - \Int_{Q_T} y \cdot \nabla{\eta} + 
	\Int_{S_N} g \eta \dst = 0,
	\quad \forall \eta \in \Ho{1}(Q_T).
	\label{eq:generalized-statement-for-v}
\end{equation}
%
%--------------------------------------------------------------------------------------%
%
The identity 
(\ref{eq:generalized-statement-for-v}) is equivalent to (\ref{eq:generalized-statement}).
Hence, $v = u$, and we see that $\Maj (v, y)$ vanishes if and 
only if
\begin{alignat}{5}
	f - v_t - \lambda v + \dvrg y & = 0  & \quad \mbox{a.e.} & \quad 
	{ (x, t)} \in Q_T, \nonumber\\
  y & = A \nabla v                     & \quad \mbox{a.e.} & \quad 
	{(x, t)} \in Q_T, \nonumber\\
	v(x, 0) & = \varphi(x)                     \; & \quad \mbox{a.e.} & \quad  
	x \in \Omega, \nonumber\\
  v & = 0                       & \quad \mbox{a.e.} & \quad 
	{(x, t)} \in S_D, \nonumber\\  
	y \cdot n & = g                      & \quad \mbox{a.e.} & \quad 
	{ (x, t)} \in S_N. 
  \label{eq:requirements-set}
\end{alignat} 
%
%--------------------------------------------------------------------------------------%
This set of requirements is fulfilled if $v$ coincides with the exact solution of the 
problem (\ref{eq:parabolic-problem-equation})--(\ref{eq:parabolic-problem-robin-boundary-condition}), 
i.e., $e = u - v = 0$ and $y$ coincides with $A \nabla u$. \hfill $\square$ 
%\end{pf}

%\end{proof}
%--------------------------------------------------------------------------------------%
%
\begin{rmk}
If $g$ is a relatively simple function, e.g., piecewise affine function, then the 
function $y$ may be selected such that $g - y \cdot n = 0$ for a.e. $(s, t) \in S_N$, 
and the constant $\Ctr$ does not appear in the estimate.
\end{rmk}
%--------------------------------------------------------------------------------------%
%
%\begin{rmk} 
%The majorant $\Maj (v, y)$ is defined if $v \in \Ho{1}(Q_T)$, $f \in \L{2}(Q_T)$, 
%$v(x, 0) \in \L{2}(\Omega)$, and $\varphi(x) \in \L{2}(\Omega)$. It is possible to 
%extend the estimate (\ref{eq:majorant-1}) to a wider set of functions. 
%\end{rmk}

\begin{rmk}
An important question, which should be discussed in the context of a posteriori
error estimation, concerns indication of local errors. 
{We note that the majorant
is presented by the sum of integrals, i.e., it automatically generates a sum of local 
quantities, which can be used as markers of local errors. In the numerical tests below, 
we show the efficiency of these error indicators. 
}
%implies such an error indicator, which can be used for adaptation of the
%space-time mesh in the process of integration. Assume, that $$ We note that the majorant can be
%represented in the form
%%
%\begin{equation}
	%\maj{} (v, y, \beta) : =
	%\| \, e(x, 0) \, \|^2_{\Omega} +
	%\frac1\delta \Int_0^T
	%\Bigg(
	%\big(1 + \beta \big) \mdI +
	%\bigg(1 + \frac{1}{\beta}\bigg) \CF^2 \mfI
	%\Bigg) \dt,
	%\label{eq:indicator-components}
%\end{equation}	
%%
%where
%%
%\begin{equation}
	%\mdI = 	\| \: y - {A} \nabla v \: \|^2_{\Omega}, \quad
	%\mfI = \| \: f - v_t - \lambda v + \dvrg y \: \|^2_{\Omega}.
	%\label{eq:indicator-components-2}
%\end{equation}
%%
%Here, $\mfI$ is `reliability term', which is necessary to provide a guaranteed upper 
%bound, and the major part of the error is usually encompassed in $\mdI$. Therefore, it 
%is natural to use the latter term as an error indicator. 
\end{rmk}
%--------------------------------------------------------------------------------------%

%--------------------------------------------------------------------------------------%
\section{Error minorant}
\label{sec:minorant}
%--------------------------------------------------------------------------------------%
Minorants of the deviation from the exact solution are well studied for elliptic 
problems, which have a variational form (see \cite{RepinDeGruyter2008, 
NeittaanmakiRepin2004} and the
literature cited therein). They provide useful information 
and allow us to judge on the quality of error majorants). Below, we derive 
computable error minorants for the evolutionary problem 
(\ref{eq:parabolic-problem-equation})--(\ref{eq:parabolic-problem-robin-boundary-condition}).

%--------------------------------------------------------------------------------------%

\begin{thm}
\label{th:theorem-minimum-of-minorant}
For any
$v, \: \eta \in \Ho{1}(Q_T)$ the following estimate holds:
\begin{equation}
	\Min (v): = 
	\sup\limits_{\eta \in \Ho{1}(Q_T)} \Bigg \{  
	\Sum_{i = 1}^{4} G_{v, i}(\eta) + F_{fg\,\varphi}(\eta) \Bigg \} \leq 
        \error_{(\underline{\nu},\, \underline{\theta},\, \underline{\zeta})} : = 
	\frac{\kappa_1}{2} \NormA{\, \nabla e } + \Bigg \|  \sqrt{ \frac{\kappa_2 + \kappa_3 \lambda}{2}} \, e \, \Bigg \|^2_{Q_T} + 
	\frac{\kappa_4}{2} \| \, e(x, T) \|^2_{\Omega}, 
	\label{eq:lower-estimate}
\end{equation}
where
\begin{alignat}{2}
G_{v, 1}(\nabla \eta) & = \Int_{Q_T} \Big( - \nabla \eta \cdot A \nabla v - \frac{1}{2 \kappa_1} |\nabla \eta|^2  \Big) \dxt, \qquad
G_{v, 2}(\eta_t) = \Int_{Q_T} \Big( \eta_t v - \frac{1}{2 \kappa_2}|\eta_t|^2 \Big) \dxt, \nonumber \\
%\end{alignat}
%
%\begin{alignat}{2}
G_{v, 3}(\eta) & = \Int_{Q_T} \lambda \Big( -  v \eta - \frac{1}{2 \kappa_3}|\eta|^2 \Big) \dxt, \qquad \; \; \;
G_{v, 4}\Big(\eta(x, T)\Big) = \Int\limits_{\Omega} \bigg( - v(x, T) \eta(x, T) - \frac{1}{2 \kappa_4}|\eta(x, T)|^2 \bigg) \dx, \nonumber \\
%\end{alignat}
%
%and
%
%\begin{equation}
F_{fg\,\varphi}(\eta) & = \Int_{Q_T} f \eta \dxt + \Int_{S_R} g \eta \dst + \Int_{\Omega} \varphi \eta(x, 0) \dx, 
%\end{equation}
\end{alignat}
parameters $\underline{\nu} = \sqrt{ \frac{\kappa_1}{2} }$, 
$\underline{\theta} = \sqrt{ \frac{\kappa_2 + \kappa_3 \lambda}{2} }$, 
$\underline{\zeta} = \sqrt{ \frac{\kappa_4}{2} }$, and 
$\kappa_1$, $\kappa_2$, $\kappa_3$, $\kappa_4 > 0$.
\end{thm}

%\begin{pf}
\noindent {\bf Proof:} 
It is not difficult to see that    
\begin{alignat}{2}
	\mathcal{M} (e) := & 
	\sup\limits_{\eta \in \Ho{1}(Q_T)} \!
	\Bigg \{
	\Int_{Q_T} \bigg( 
	\nabla \eta \cdot A \nabla e - \frac{1}{2 \kappa_1} |\nabla \eta|^2 - 
        \eta_t e - \frac{1}{2 \kappa_2} |\eta_t|^2 + 
        \lambda \Big( e \eta - \frac{1}{2 \kappa_3} |\eta|^2 \Big) 
	\bigg) \dxt + \nonumber\\
	%\qquad \qquad \qquad \qquad \qquad \qquad \qquad \qquad \qquad \qquad \qquad \qquad
	& \qquad  \qquad \qquad \qquad  \qquad \qquad \qquad \qquad \qquad \qquad
	\Int_\Omega \Big( e(x, T) \eta (x, T) - \frac{1}{2 \kappa_4} |\eta (x, T)|^2 \Big)\dx \Bigg \} \leq \nonumber \\ 	
	\sup\limits_{\eta \in \Ho{1}(Q_T)} \Bigg \{ \Int\limits_{Q_T} \bigg( & 
	\nabla \eta \cdot A \nabla e - \frac{1}{2 \kappa_1} |\nabla \eta|^2 \bigg) \dxt \Bigg \} + 
	\sup\limits_{\eta \in \Ho{1}(Q_T)} \Bigg \{ \Int\limits_{Q_T} \Big( - \eta_t e - \frac{1}{2 \kappa_2} |\eta_t|^2 \Big) \dxt \Bigg \} + \nonumber \\
	& \sup\limits_{\eta \in \Ho{1}(Q_T)} \Bigg \{ \Int\limits_{Q_T} \lambda \Big( e \eta - \frac{1}{2 \kappa_3} |\eta|^2 \Big) \dxt \Bigg \} + 
	\sup\limits_{\eta (x, T) \in \Ho{1}(\Omega)} \Bigg \{ \Int\limits_{\Omega} \Big( e(x, T) \eta (x, T) - \frac{1}{2 \kappa_4} |\eta|^2 \Big)\dx \Bigg \}. \nonumber
	%\label{eq:quadratic-func-inequality}the 
\end{alignat}
Since
\begin{alignat}{2}
\sup\limits_{\eta \in \Ho{1}(Q_T)} \Bigg \{ \Int\limits_{Q_T} \bigg( \nabla \eta \cdot A \nabla e - \frac{1}{2 \kappa_1} |\nabla \eta|^2 \bigg) \dxt \Bigg \} & \leq \frac{\kappa_1}{2} \NormA{\nabla e } \,, \nonumber \\
\sup\limits_{\eta \in \Ho{1}(Q_T)} \Bigg \{ \Int\limits_{Q_T} \Big( - \eta_t e - \frac{1}{2 \kappa_2} |\eta_t|^2 \Big) \dxt \Bigg \} & \leq \frac{\kappa_2}{2} \| e \|^2_{Q_T}, \nonumber \\
\sup\limits_{\eta \in \Ho{1}(Q_T)} \Bigg \{ \Int\limits_{Q_T} \lambda \Big( e \eta - \frac{1}{2 \kappa_3} |\eta|^2 \Big) \dxt \Bigg \} & \leq \frac{\kappa_3}{2} \big\| \sqrt{\lambda} e \big\|^2_{Q_T}, \nonumber \\
\sup\limits_{\eta (x, T) \in \Ho{1}(\Omega)} \Bigg \{ \Int\limits_{\Omega} \Big( e(x, T) \eta (x, T) - \frac{1}{2 \kappa_4} |\eta (x, T)|^2 \Big)\dx \Bigg \} & \leq \frac{\kappa_4}{2} \| e(x, T) \|^2_{\Omega}, \nonumber
\end{alignat}
we find that 
\begin{equation}
	%\sup\limits_{\eta \in \Ho{1}(Q_T)}
	%\Bigg \{
	%\Int\limits_{Q_T} \bigg(
	%\nabla \eta \cdot A \nabla e - \frac{1}{2 \kappa_1} |\nabla \eta|^2 
  %- \eta_t e - \frac{1}{2 \kappa_2} |\eta_t|^2 + 
	%\lambda \Big( e \eta - \frac{1}{2 \kappa_3} |\eta|^2 \Big) 
	%\bigg) \dxt + \\
	%\Int\limits_{\Omega} \Big( e(x, T) \eta (x, T) - \frac{1}{2 \kappa_4} |\eta|^2 \Big)\dx \Bigg \} 
	\mathcal{M} (e) 
	%[e]_{(\nu, \, \theta, \, \zeta)} : = 
	\leq \error_{(\underline{\nu},\, \underline{\theta},\, \underline{\zeta})},
	%: = 
	%\frac{\kappa_1}{2} \, \NormA{\nabla e } + \, \Bigg \| \,  \sqrt{ \frac{\kappa_2 + \kappa_3 \lambda}{2}}  \, e  \, \Bigg  \|^2_{Q_T} + 
	%\frac{\kappa_4}{2} \, \| \, e(x, T) \, \|^2_{\Omega}.
	\label{eq:quadratic-func-inequality-1}
\end{equation}
where $\error_{(\underline{\nu},\, \underline{\theta},\, \underline{\zeta})}$ is defined
in (\ref{eq:lower-estimate}).
%--------------------------------------------------------------------------------------%
%
On the other hand, by using (\ref{eq:generalized-statement}), we see that for any 
$\eta$ {the} functional 
\begin{equation}
	%\sup\limits_{\eta \in \Ho{1}(Q_T)} 
	%\Bigg \{
	%\Int_{Q_T} 
  %\bigg(
	%\nabla \eta \cdot A \nabla e - \frac{1}{2 \kappa_1} |\nabla \eta|^2 
	%- \eta_t e - \frac{1}{2 \kappa_2} |\eta_t|^2 + 
	%\lambda \Big( e \eta - \frac{1}{2 \kappa_3} |\eta|^2 \Big) \bigg) \dxt + 
	%\IntO \Big( e(x, T) \eta (x, T) - \frac{1}{2 \kappa_4} |\eta|^2 \Big)\dx \Bigg \} 
	\mathcal{M} (e)
	= \\
	\sup\limits_{\eta \in \Ho{1} (Q_T)} \Bigg \{  
	\Sum_{i = 1}^{4} G_{v, i} + F_{fg\,\varphi}(\eta) \Bigg \}
\end{equation}
generates 
{a} lower bound of the norm 
$\error_{(\underline{\nu},\, \underline{\theta},\, \underline{\zeta})}$.
Thus, we arrive at (\ref{eq:lower-estimate}). \hfill $\square$
%
%\end{pf}

%--------------------------------------------------------------------------------------%
\section{Incremental form of the estimates}
\label{sec:incremental}
%--------------------------------------------------------------------------------------%

%--------------------------------------------------------------------------------------%
Let $\mathcal{T}_K = \bigcup\limits_{k = 0}^{K-1}[t^k, t^{k + 1}]$ be a mesh selected 
on $[0,\; T]$, so that $Q_T$ can be represented in the form
\begin{equation}
\overline{Q}_T = \bigcup\limits_{k = 0}^{K-1} \overline{Q^k}, \quad 
       Q^k := ( t^k, t^{k+1} ) \times \Omega,
\end{equation}
and $\mathcal{T}_{N_1 \times .. \times N_d}$ be a mesh selected on $\Omega$. Therefore,  
$\Theta_{K \times N_1 \times ... \times N_d} = 
\mathcal{T}_K \times \mathcal{T}_{N_1 \times ... \times N_d}$
denotes the mesh on $Q_T$.

Computational methods developed for parabolic type of problems often use incremental 
(semi-discrete) schemes. Below, we discuss forms of the majorants and minorants 
adapted to this class of methods. They lead to estimates, which evaluate errors on 
each interval $[t^k, \;t^{k+1}]$ and accumulate the overall error. 
In this section, we assume that
%
%In computations, it is convenient  to have incremental form of the estimates, in 
%which errors are evaluated for each time-step $[t^k, t^{k+1}]$. Below, we discuss 
%such forms assuming that 
%
$A = I$, $\lambda = \lambda(x)$, and $S_T = S_D$ (these assumptions 
are made in order to simplify the notation). 

We set  
\begin{equation}
	\alpha_1 \big( \beta(t), \delta \big) = \frac1\delta \Bigg( 1 + \frac{1}{\beta(t)} \Bigg), \quad 
	\alpha_2 \big( \beta(t), \delta \big) = \frac1\delta \Big( 1 + \beta(t) \Big),
\end{equation}
where $\beta(t) \in L^{\infty}_{\hat{\nu}} (t^k, t^{k+1}) := 
\Bigg\{ \: \beta(t) \in L^\infty_{\hat{\nu}} (t^k, t^{k+1})\; \big|
\;\beta(t) \geq \hat{\nu} > 0\; \mbox{for a.e. } t \in (t^k, t^{k+1})\: \Bigg\}$.

Assume that $v$ is computed by a simple semi-discrete approximation method on interval
$[t^k, t^{k+1}]$ of length $\tau = t^{k+1} - t^k$ (see, e.g.,
 \cite{Thomee2006, Braess2001, Johnson2009, LeVeque2007}):
\begin{equation*}
v = v^{k} \: \frac{t^{k+1} - t}{\tau} + v^{k+1} \: \frac{t - t^{k}}{\tau}, \quad v_t = \frac{v^{k+1} - v^{k}}{\tau}, \quad
y = y^{k} \: \frac{t^{k+1} - t}{\tau} + y^{k+1} \: \frac{t - t^{k}}{\tau}, \quad f = f^{k} \: \frac{t^{k+1} - t}{\tau} + f^{k+1} \: \frac{t - t^{k}}{\tau}.
%\label{eq:v-y-linerization-of-time}
\end{equation*}
Next, we define 
\begin{alignat}{2}
	\R^k_d(x) & = \R^k_d = { y^k - {\nabla v}^k}, \nonumber\\
	\R^k_f(x) & = \R^k_f = \dvrg y^k 
	{- \lambda^k v^k } 
	+ f^k - \frac{ v^{k + 1} - v^{k}}{\tau}, \nonumber \\
	\R^{k+1}_f(x) & = \R^{k+1}_f = \dvrg y^{k+1} 
	{	- \lambda^{k+1} v^{k+1}} + f^{k+1} - \frac{  v^{k + 1} - v^{k}}{\tau},
\end{alignat}
which yields
\begin{equation}
\R_d(x, t) = \R^k_d \: \frac{t^{k+1} - t}{\tau} + 
             \R^{k+1}_d \: \frac{t - t^{k}}{\tau}, \quad 
\R_f(x, t) = \R^k_f \: \frac{t^{k+1} - t}{\tau} + 
             \R^{k+1}_f \: \frac{t - t^{k}}{\tau}.
\end{equation}
{Since
\begin{alignat*}{2}
\IntT & (t - t^k) \dt = \IntT (t^{k+1} - t) \dt = \frac{\tau^2}{2}, \quad 
\frac12 \IntT (t - t_k)^2 \dt = \frac12 \IntT (t_{k+1} - t)^2 \dt = 
\IntT (t - t_k)(t_{k+1} - t) = \frac{\tau^3}{6}, \\
\IntT & (t - t_k)^2 (t_{k+1} - t) \dt = 
\IntT (t_{k+1} - t)^2(t - t_k) \dt = \frac{\tau^4}{12}, \quad
\IntT (t - t_k)^2 (t_{k+1} - t)^2\dt = \frac{\tau^5}{30}, \\
%\end{alignat}
%
%\begin{alignat}{2}
\IntT & \frac{(t_k + t_{k+1} - 2t)(t - t_k)}{\tau} \dt = - \frac{\tau^2}{6}, \quad
\IntT \frac{(t_k + t_{k+1} - 2t)(t_{k+1} - t)}{\tau} \dt = \frac{\tau^2}{6}, \quad \\
\IntT & (t_k + t_{k+1} - 2t) \dt = 0, \quad 
\IntT (t_k + t_{k+1} - 2t)^2 \dt = \frac{\tau^3}{3}, 
\end{alignat*}
we find that for $\overline{Q^0}$
} 
\begin{multline}
	\majk{(0)} = 
	\| \: e(x, t^0) \: \|^2_{\Omega} + 
	\Big(1 + \beta^0 \Big)  \IntO \Int_{t^0}^{t^1}   
	\Bigg( \R^0_d \frac{t^1 - t}{\tau} + \R^{1}_d \frac{t - t^0}{\tau} \Bigg)^2 \dt \dx + \\
	\qquad \qquad \qquad \qquad \qquad \qquad \qquad \qquad \qquad \qquad
	\CF^2 \Bigg( 1 + \frac{1}{\beta^0} \Bigg) \IntO \Int_{t^0}^{t^1} 
	\Bigg( \R^0_f \frac{t^1 - t}{\tau} + \R^1_f \frac{t - t^0}{\tau} \Bigg)^2 \dt \dx  = \\
	\| \: \varphi - v^0 \: \|^2_{\Omega} +
	\frac{\tau}{3} \Bigg \{ 
	\Big(1 + \beta^0 \Big) \IntO \Bigg( \Big(\R^0_d\Big)^2 + \R^0_d \R^1_d + \Big(\R^1_d \Big)^2 \Bigg) \dx + 
	\qquad \qquad \qquad \qquad \qquad \qquad\\
	\CF^2 \Bigg( 1 + \frac{1}{\beta^0} \Bigg) \IntO \Bigg(  \Big(\R^0_f\Big)^2 + \R^0_f \R^1_f + \Big(\R^1_f\Big)^2 \Bigg) \dx 
	\Bigg \}.
	\label{eq:error-majorant-time-level-0}
\end{multline} 
%
%--------------------------------------------------------------------------------------%
%
If $\majk{(k-1)}$ is the majorant related to 
$\overline{Q}_{t^k} = \bigcup\limits_{j = 0}^{k-1} \overline{Q^j}$, 
then $\majk{(k)}$ on $\overline{Q}_{t^{k+1}}$ can be computed by the recurrent 
formula
\begin{multline}
	\majk{(k)}  = \: \majk{(k-1)} + 
	\frac{\tau}{3} \Bigg \{ 
	\Big(1 + \beta^k \Big) \IntO
	\Bigg( \Big(\R^k_d\Big)^2 + \R^k_d \R^{k+1}_d + \Big(\R^{k+1}_d\Big)^2 \Bigg) \dx + \\
	\CF^2 \Bigg( 1 + \frac{1}{\beta^k} \Bigg) \IntO 
	\Bigg( \Big(\R^k_f\Big)^2 + \R^k_f \R^{k+1}_f + \Big(\R^{k+1}_f\Big)^2 \Bigg) \dx 
	\Bigg \}.
	\label{eq:error-majorant-time-level-k}
\end{multline} 
%
%--------------------------------------------------------------------------------------%
%
Analogously, we deduce similar relation for the minorant.
We set 
\begin{alignat}{2}
\eta & = \eta^{k} \: \frac{t^{k+1} - t}{\tau} + \eta^{k+1} \: \frac{t - t^k}{\tau} + \alpha (t - t^k) (t^{k+1} - t), \quad 
\eta_t = \frac{\eta^{k+1} - \eta^k}{\tau} + \alpha (t^k + t^{k+1} - 2t), \nonumber \\
\nabla \eta & = \nabla \eta^{k} \: \frac{t^{k+1} - t}{\tau} + \nabla \eta^{k+1} \: \frac{t - t^{k}}{\tau} + \nabla \alpha (t - t^k) (t^{k+1} - t).
\label{eq:eta-approximation}
\end{alignat}
Then, on $\overline{Q}_{t^{k+1}}$ the minorant is presented as
\begin{equation}
	\mink{(k)} : = \mink{(k-1)} + 
	\Sum_{i = 1}^{4} G^k_{v, i}(\eta) + F^k_{fg\, \varphi}(\eta), 
	\label{eq:incremental-lower-estimate}
\end{equation}
where $\mink{(k-1)}$ is constructed on $\overline{Q}_{t^k}$ and
\begin{alignat}{2}
G^k_{v, 1} & = - \frac{\tau}{3} 
                 \IntO	\Bigg (
								 \nabla \eta^{k} \nabla v^{k} + \nabla \eta^{k+1} \nabla v^{k+1} + 
								        \frac12 \bigg( \nabla \eta^{k+1} \nabla v^{k} + \nabla \eta^{k} \nabla v^{k+1} \bigg) +
                        \frac{\tau^2}{4} \nabla \alpha \bigg(\nabla v^{k} + \nabla v^{k+1} \bigg) + \nonumber \\
					 & \qquad \qquad \qquad \qquad \frac{1}{2 \kappa_1} 
						 \bigg( \Big(\nabla \eta^{k}\Big)^2 + \Big(\nabla \eta^{k+1} \Big)^2 + \nabla \eta^{k} \nabla \eta^{k+1} + 
						 \frac{\tau^2}{2} \nabla \alpha \Big(\nabla \eta^{k} + \nabla \eta^{k+1}\Big)
								    + \frac{\tau^4}{10} \nabla \alpha^2 \bigg) \Bigg ) \dx, \\
G^k_{v, 2} & = { \frac12 \IntO \Bigg( 
														 \Big( v^{k} + v^{k+1} \Big) \Big(\eta^{k+1} - \eta^{k}\Big) + 
										                \alpha^k \frac{\tau^2}{3} \Big(v^{k} - v^{k+1} \Big) 
															- \frac{1}{\kappa_2} 
																\Bigg( \frac{ \big(\eta^{k+1} - \eta^{k}\big)^2}{\tau} + 
																       \Big(\alpha^k\Big)^2 \frac{\tau^3}{3} \Bigg) 
															\Bigg) \dx,} \\ 																				
G^k_{v, 3} & = - \frac{\tau}{3} 
               \IntO {\lambda} \Bigg( \eta^{k} v^{k} + \eta^{k+1} v^{k+1} + \frac12 \Big(\eta^{k+1} v^{k} + \eta^{k} v^{k+1} \Big) 
                             + \frac{\tau^2}{4} \alpha \Big(v^{k} + v^{k+1}\Big)  + \nonumber \\
					 & \qquad \qquad \qquad \qquad \qquad \qquad \quad 
					   \frac{1}{2 {\kappa_3}} 
						 \bigg( \Big( \eta^{k} \Big)^2 + \Big( \eta^{k+1} \Big)^2 + \eta^{k} \eta^{k+1} + \frac{\tau^2}{2} \alpha \Big( \eta^{k} + \eta^{k+1} \Big)
								    + \frac{\tau^4}{10} \alpha^2 \bigg)	\Bigg)\dx, \\
G^k_{v, 4} & = \IntO	\Bigg( - \eta^{k+1} v^{k+1} - \frac{1}{2 \kappa_4} \Big(\eta^{k+1}\Big)^2 \Bigg) \dx, 
\end{alignat}

\begin{alignat}{2}
F^k_{fg\,\varphi}   & = {\frac{\tau}{3}} \IntO \bigg( \eta^{k} f^{k} + \eta^{k+1} f^{k+1} + \frac12 \Big( \eta^{k+1} f^{k} + \eta^{k} f^{k+1} \Big) + \frac{\tau^2}{4} \alpha \Big(f^{k} + f^{k+1} \Big) \bigg) \dx + \nonumber \\
           & \qquad \qquad \qquad \qquad \qquad \qquad 
					  {\frac{\tau}{3}} 
						\IntGammaN \bigg( \eta^{k} g^{k} + \eta^{k+1} g^{k+1} + \frac12 \Big( \eta^{k+1} g^{k} + \eta^{k} g^{k+1} \Big) + \frac{\tau^2}{4} \alpha \Big(g^{k} + g^{k+1} \Big) \bigg)\ds.
\label{eq:g-f}
\end{alignat}
\begin{rmk}

{We note that the presented incremental forms of the estimates 
(\ref{eq:error-majorant-time-level-0})--(\ref{eq:error-majorant-time-level-k}) and 
(\ref{eq:incremental-lower-estimate})--(\ref{eq:g-f})
are valid only for the first order time discretization scheme. 
The corresponding higher order schemes can be derived by similar arguments, 
but the estimates will have a more complicated form. However, this subject is beyond 
the framework of this paper and will be considered in a subsequent
publications. 
%For higher order time discretization schemes, these estimates will have to be %accordingly modified.
}
Also, we note that in the case of an oscillating right-hand side, we can combine the current estimate with the majorant or minorant of the modeling error (see, e.g., 
\cite{RepinSamrowskiSauter2012}).
\end{rmk}

%--------------------------------------------------------------------------------------%
\section{Numerical tests}
\label{sec:numerics}
%--------------------------------------------------------------------------------------%

%--------------------------------------------------------------------------------------%
% example-3
%--------------------------------------------------------------------------------------%

\begin{example}
\label{ex:example-3}
\rm
We begin with a relatively simple problem
where $\Omega=(0, 1)$, $T = 10$, { $\partial \Omega = \Gamma_D$}, 
$u=0$ on
$S_D$, $A = I$, $\varphi = x\, \left(1 - x \right)$, { $\lambda(x) = 0$,}
and $f = 2\, t ( 1 + \,t) - x\, \left(2\, t + 1\right)\, \left(x - 1\right) + 2$.
The corresponding exact solution is \linebreak
$u = x\, \left(1 - x \right)\, \left(t^2 + t + 1\right)$. 

%
%--------------------------------------------------------------------------------------%
%
%In order to have a realistic presentation of the accuracy, we normalize absolute values
%of the error by the energy norm of the exact solution $[u]^2$.
The quality of the error estimates is measured by efficiency indexes:
%
%\begin{alignat}{2}
%\Ieff \left(\maj{} \right) 			       & := \frac{\overline{\mathrm M}}{[e]_{(\nu,\, \theta,\, \zeta)}} \geq 1,  \\
%\Ieff \left(\Min \right) 				       & := \frac{\underline{\mathrm M}}{[e]_{(\nu,\, \theta,\, \zeta)}} \leq 1,  \\
%\Ieff \left(\frac{\maj{}}{\Min}\right) & := \frac{\overline{\mathrm M}}{\underline{\mathrm M}}\geq 1.
%\end{alignat}
{
\begin{equation}
\Ieffmaj := \frac{\overline{\mathrm M}}{[e]_{(\nu,\, \theta,\, \zeta)}} \geq 1,  \quad
\Ieffmin := \frac{\underline{\mathrm M}}{[e]_{(\nu,\, \theta,\, \zeta)}} \leq 1,  \quad	
\Ieff := \frac{\overline{\mathrm M}}{\underline{\mathrm M}} \geq \Ieffmaj \geq 1 \geq \Ieffmin.
\end{equation}
}
In order to have a realistic presentation of the accuracy, we normalize
values of the error and estimates by the energy norm of the exact solution $[u]^2$,
i.e.,  we compare the relative values
$\frac{\maj{}}{[u]^2}$ and $\frac{\Min}{[u]^2}$
with the relative value of the true error
$\frac{\error_{(\nu,\, \theta,\, \zeta)}}{[u]^2}$.
Table \ref{tab:example-3-error-maj-ieff-for-different-delta} shows these quantities
for
%the relative 
%error and respective relative majorant associated with 
different values of the parameter 
$\delta$ (approximate solution was computed on the uniform mesh
$\Theta_{\:K \times N_1}$ = $\Theta_{\:40 \times 40}$). 
{
We see that within the interval $[0.5, \,1.5]$ the efficiency of the majorant is not
very sensitive with respect to the parameter $\delta$. Results of other examples suggest
 similar conclusions. Therefore, in subsequent tests we set $\delta$ = 1
%In further examples, we set 
%$\delta = 1$ 
in order to obtain the optimal value of the majorant.
}
\begin{table}[!ht]
	\centering
	\footnotesize
	\begin{tabular}{cccc}
	\midrule
	$\delta$ & $\error_{\left(\sqrt{2-\delta}, \, 0, \, 1 \right)} / [u]^2$
					 & $\maj{} / [u]^2$
					 & $\Ieffmaj$ \\
	\midrule
	0.5 & 6.37e-004 & 8.66e-004 & 1.17 \\
		1 & 6.28e-004 & 6.39e-004 & 1.01 \\
	1.5 & 6.01e-004 & 8.16e-004 & 1.17 \\
	\midrule
	\end{tabular}
	\caption{{Example \ref{ex:example-3}}.
	The relative error, majorant, and its efficiency index with respect to parameter
	$\delta$ 
	{for $t = T$}.}
	\label{tab:example-3-error-maj-ieff-for-different-delta}
\end{table}
%--------------------------------------------------------------------------------------%

Growth of the error $\log \, \error_{(1, \, 0, \, 1)}$ and the
majorant $\log \, \maj{}$ in logarithmic scale is depicted on 
Fig. \ref{fig:example-3-log-abs-e-maj}. We see 
that the majorant reproduces the error quite accurately. Table 
\ref{tab:example-3-majorant-terms} presents different components of the majorant.
The term
${{\| \: y - \nabla v \: \|}^2_{\overline{Q}_{{t\,}^k}}}$ contains the main part of the majorant
and represents the error term
${{\mid\mid\mid\! \nabla(u - v) \!\mid\mid\mid}^2_{ \overline{Q}_{{t\,}^k} }}$
quite accurately so that the efficiency index is close to $1$ for any $t \in [0, T]$. 
%Since the
%term ${{\| \: f - v_t + \dvrg y \:\|}^2_{ \overline{Q}_{{t\,}^k} }}$ is relatively small, 
%the obtained efficiency index is close to one.

%
\begin{figure}[!ht]
	\centering
  \includegraphics[scale=0.9]{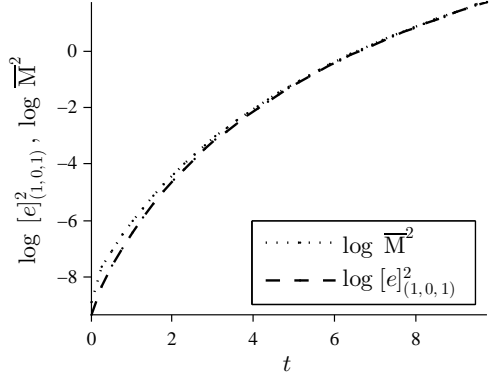}
  \caption{Example \ref{ex:example-3}.
	The error and the majorant with respect to time. }
	\label{fig:example-3-log-abs-e-maj}
\end{figure}
%
%--------------------------------------------------------------------------------------%
%
\begin{table}[!th]
	\centering
	\footnotesize
	\begin{tabular}{cccccc}
	\midrule
	$t^k$ & ${{\mid\mid\mid \! \nabla(u - v) \! \mid\mid\mid}^2_{ \overline{Q}_{{t \,}^k} }} $
			& $\| \: (u - v)(x, t^k) \:\|^2_{\Omega}$
			& ${{\| \: y - \nabla v \:\|}^2_{ \overline{Q}_{{t \,}^k}}}$
			& ${{\| \: f - v_t + \dvrg y \:\|}^2_{ \overline{Q}_{{t \,}^k} }}$ 
			& $\Ieffmaj$\\
	\midrule
	1.03 & 4.28e-04 & 1.52e-07 & 4.24e-04 & 2.35e-04 & 1.28\\
	2.05 & 2.49e-03 & 8.44e-07 & 2.49e-03 & 1.97e-04 & 1.13\\
	3.08 & 9.02e-03 & 2.94e-06 & 9.02e-03 & 1.94e-04 & 1.07\\
	4.10 & 2.43e-02 & 7.70e-06 & 2.43e-02 & 1.94e-04 & 1.05\\
	5.13 & 5.41e-02 & 1.68e-05 & 5.41e-02 & 2.85e-04 & 1.03\\
	6.15 & 1.06e-01 & 3.25e-05 & 1.06e-01 & 2.85e-04 & 1.02\\
	7.18 & 1.88e-01 & 5.71e-05 & 1.88e-01 & 2.86e-04 & 1.01\\
	8.21 & 3.10e-01 & 9.38e-05 & 3.10e-01 & 1.93e-04 & 1.01\\
	9.23 & 4.86e-01 & 1.46e-04 & 4.86e-01 & 3.04e-04 & 1.01\\
	10.00 & 6.60e-01 & 1.97e-04 & 6.60e-01 & 2.35e-04 & 1.01\\
	\midrule
	\end{tabular}
	\caption{{Example \ref{ex:example-3}}.
	Two terms of the error, two terms of the majorant, and corresponding efficiency index 	
	with respect to time.}
	\label{tab:example-3-majorant-terms}
\end{table}
%--------------------------------------------------------------------------------------%

The normalized quantities  $\frac{\error_{(1, \: 0, \: 1)}}{[u]^2}$ and
$\frac{\maj{}}{[u]^2}$ (as functions of time) are depicted in Fig. 
\ref{fig:example-3-rel-e-maj}. In Fig. \ref{fig:example-3-rel-e-min}, we show the 
value
$\frac{\error_{(1/\sqrt{2}, \: 1/\sqrt{2}, \: 1/\sqrt{2})}}{[u]^2}$ and the 
corresponding 
minorant $\frac{\Min{}}{[u]^2}$. We see that the computed two-sided bounds of the 
error are efficient and guaranteed for any $t^k \in [0, T]$.

\begin{figure}[!ht]
	\centering
  \subfloat[]{
	\includegraphics[scale=0.9]{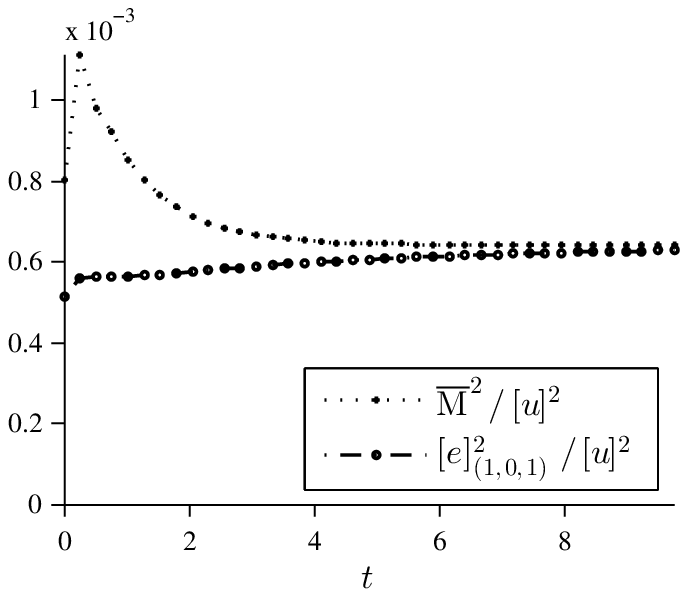}
	\label{fig:example-3-rel-e-maj}
	}
	\subfloat[]{
	\includegraphics[scale=0.9]{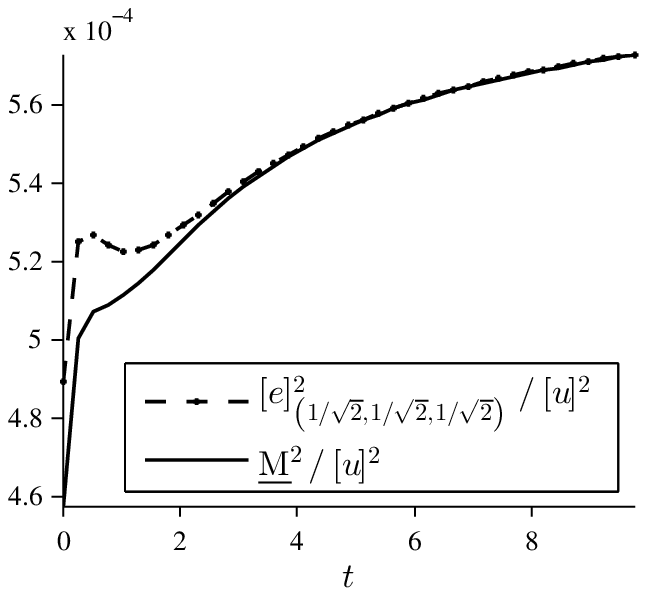}
	\label{fig:example-3-rel-e-min}
	}
  \caption{Example \ref{ex:example-3}.
	(a) The relative error and majorant, (b) the relative error and minorant with respect to time.
	}
\end{figure}

If the weights of the error terms are selected in the special way, then the majorant and the minorant
generate two-sided bounds for the same error norm (\ref{eq:energy-norm-for-reaction-diff-evolutionary-problem}).
For any $\varkappa \in (0, 2 - \delta)$, we have
\begin{multline}
(2 - \delta)\NormA{\nabla e} + \| \: e(x, T) \: \|^2_{\Omega} =
(2 - \delta - \varkappa) \, \NormA{\nabla e} + \varkappa \, \NormA{\nabla e} + \| \: e(x, T) \: \|^2_{\Omega} \geq \\
(2 - \delta - \varkappa) \, \NormA{\nabla e} + \frac{\varkappa \nu_1 }{\CF^2} \, \| \: e \: \|^2_{Q_T} + \| \: e(x, T) \: \|^2_{\Omega}.
\label{eq:majorant-norm}
\end{multline}
%
%To balance right-hand side of (\ref{eq:majorant-norm}) with the norm
%
%\begin{equation}
%\frac{\kappa_1}{2} \NormA{\nabla e } + \frac{\kappa_2}{2} \| e \|^2_{Q_T} +
%\frac{\kappa_4}{2} \| e(x, T) \|^2_{\Omega},
%\label{eq:minorant-norm}
%\end{equation}
%
%we have to select the parameter as follows:
%
%\begin{equation}
%\frac{\kappa_1}{2} = 2 - \delta - \varkappa, \quad
%\frac{\kappa_2}{2} = \frac{ \varkappa \nu_1^2 }{\CF^2}, \quad \mbox{and} \quad
%\frac{\kappa_4}{2} = 1.
%\end{equation}
%
Then, the left-hand side is estimated from above by the majorant $\maj{}$ 
(see Theorem \ref{eq:majorant-1})
and the right-hand side is estimated from below by the minorant $\Min$ (see Theorem 
\ref{eq:lower-estimate}) with the corresponding weights.
Henceforth, 
we set $\delta = 1$, then weighted error norms are 
denoted by $\error_{(\tilde{\nu},\, \tilde{\theta},\, \tilde{\zeta})}$, where 
$\tilde{\nu} = \sqrt{1 - \varkappa}$,
$\tilde{\theta} = \frac{\sqrt{ \varkappa}}{\CF}$, and
$\tilde{\zeta} = 1$. In minorant, 
$\kappa_1 = 2 \, (1 - \varkappa)$, 
$\kappa_2 = \frac{2 \varkappa}{\CF^2}$, and $\kappa_4 = 2$. By changing $\varkappa$, 
we obtain different weighted norms of the error, 
which have computable two-sided error bounds. 
In Table \ref{tab:balanced-error-i-effs-from-varkappa}, we present the efficiency index of 
these two-sided bounds for different $\varkappa$. In Fig. \ref{fig:example-3-rel-e-min-maj-from-ae}, we 
depict the behavior of 
$\frac{\error_{(\tilde{\nu},\, \tilde{\theta},\, \tilde{\zeta})}}{[u]^2}$,
$\frac{\Min{}}{[u]^2}$, and $\frac{\maj{}}{[u]^2}$ with respect to 
time for $\varkappa$ = $5 \cdot 10^{\minus 1}$ \linebreak
and $\varkappa$ = $10^{\minus 3}$. 
%from which it is 
%obvious that the estimates provide more accurate results if $\varkappa$ goes to zero.

%--------------------------------------------------------------------------------------%

\begin{table}[!t]
\footnotesize
\centering
\begin{tabular}{ccccccc}
\midrule
$\:$ & \multicolumn{2}{c}{$\varkappa = 5 \cdot 10^{ \minus 1}$}
     & \multicolumn{2}{c}{$\varkappa = 10^{ \minus 2}$}
		 & \multicolumn{2}{c}{$\varkappa = 10^{ \minus 3}$}\\
\cmidrule{2-7}
$t^k$ & $\error_{(\tilde{\nu},\, \tilde{\theta},\, \tilde{\zeta})} / [u]^2 $ & $ \Ieff = \frac{\overline{\mathrm M}}{\underline{\mathrm M}}$ 
% $ \Ieff \left(\frac{\maj{}}{\Min}\right) $
& $\error_{(\tilde{\nu},\, \tilde{\theta},\, \tilde{\zeta})} / [u]^2 $ 
& $\Ieff = \frac{\overline{\mathrm M}}{\underline{\mathrm M}}$
% $ \Ieff \left(\frac{\maj{}}{\Min}\right) $
& $\error_{(\tilde{\nu},\, \tilde{\theta},\, \tilde{\zeta})} / [u]^2 $ 
& $\Ieff = \frac{\overline{\mathrm M}}{\underline{\mathrm M}}$ \\
% $ \Ieff \left(\frac{\maj{}}{\Min}\right) $ 
% \\
% $ \Ieff \left(\frac{\maj{}}{\Min}\right) $ 
\midrule
0.26 & 2.95e-04 & 3.02  & 5.10e-04 & 1.69 & 5.14e-04 & 1.69 \\
1.03 & 2.98e-04 & 3.36  & 5.58e-04 & 1.69 & 5.63e-04 & 1.68 \\
2.05 & 2.91e-04 & 2.60  & 5.66e-04 & 1.31 & 5.71e-04 & 1.30 \\
3.08 & 2.95e-04 & 2.31  & 5.79e-04 & 1.17 & 5.85e-04 & 1.16 \\
4.10 & 3.01e-04 & 2.19  & 5.91e-04 & 1.11 & 5.96e-04 & 1.10 \\
5.13 & 3.05e-04 & 2.13  & 6.00e-04 & 1.08 & 6.05e-04 & 1.07 \\
6.15 & 3.08e-04 & 2.10  & 6.07e-04 & 1.06 & 6.12e-04 & 1.05 \\
7.18 & 3.11e-04 & 2.08  & 6.12e-04 & 1.05 & 6.18e-04 & 1.04 \\
8.21 & 3.13e-04 & 2.06  & 6.16e-04 & 1.04 & 6.22e-04 & 1.03 \\
9.23 & 3.15e-04 & 2.05  & 6.20e-04 & 1.04 & 6.25e-04 & 1.03 \\
10.00 & 3.16e-04 & 2.04  & 6.22e-04 & 1.03 & 6.27e-04 & 1.02 \\
\midrule
\end{tabular}
\caption{Example \ref{ex:example-3}.
				 The error and efficiency index 
				 $\Ieff = \frac{\overline{\mathrm M}}{\underline{\mathrm M}}$ with respect 
				 to time for different parameter $\varkappa$.}
\label{tab:balanced-error-i-effs-from-varkappa}
\end{table}
\begin{figure}[!t]
	\centering
	\subfloat[]{
	\includegraphics[scale=0.9]{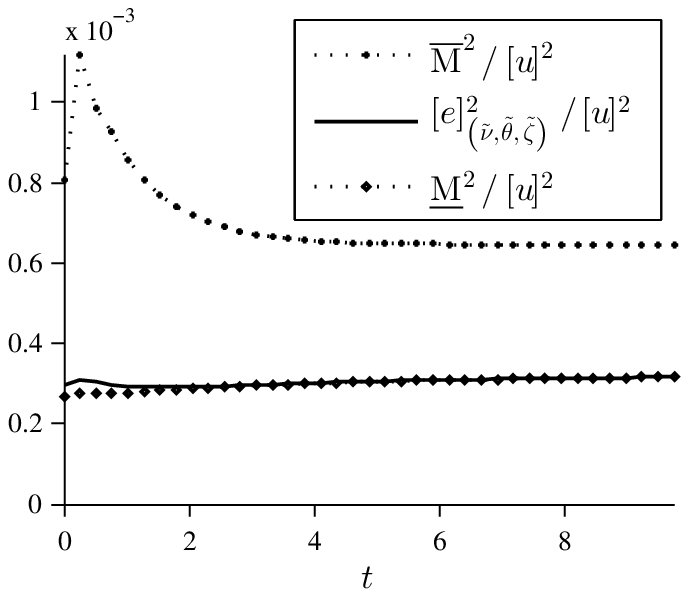}
	\label{fig:example-3-rel-e-maj-from-ae}
	}
	\subfloat[]{
	\includegraphics[scale=0.9]{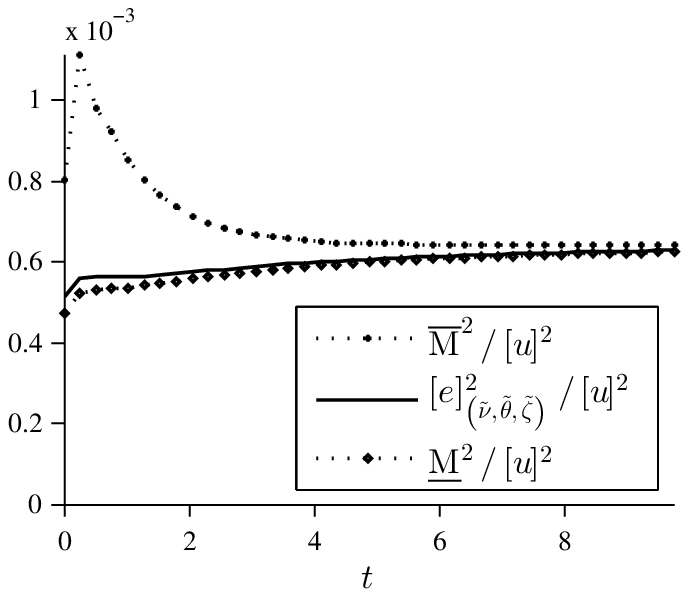}
	\label{fig:example-3-rel-e-min-from-ae}
	}
  \caption{Example \ref{ex:example-3}.
	Two-sided bounds of the error with respect to time for (a) 
	$\varkappa = 10^{ - 1}$  and (b) $\varkappa = 5 \cdot 10^{ - 3}$.}
  \label{fig:example-3-rel-e-min-maj-from-ae}
\end{figure}

%--------------------------------------------------------------------------------------%

The efficiency of the majorant (minorant) depends on the selection of
$y \in Y_{\dvrg}(Q_T)$ $ \left( \eta \in \Ho{1}(Q_T) \right)$. For elliptic problems,
this question is well studied, and there exist methods which are able to reconstruct $y$
using the computed flux $A \nabla v$ (see, e.g., \cite{RepinTomar2011, Valdman2009,
KleissTomar2012}). We follow the same technique and obtain a suitable $y$ ($\eta$)
by local minimization of $\maj{}$ (maximization of $\Min$). The corresponding results 
are collected in Table \ref{tab:example-3-maj-I-minimization-with-y}. The first two 
columns show the 
performance of $\maj{}$ with the most coarse reconstruction, in which $y$ is obtained 
by a simple patch-wise averaging of the numerical flux $A \nabla v$. This procedure 
is very cheap and does not lead to noticeable computational costs. The columns 5 and 
6 show the estimates obtained after minimization process, in which $y$ was computed 
by means of a patch-wise minimization of $\maj{}$. The last two 
columns present the results related to the best possible reconstruction of the flux 
obtained by adding patch-wise based bubble functions.
%This procedure requires additional 
%time and provides practically the same results as in column 7 and 8. 
Table \ref{tab:example-3-min-maximization-with-eta} presents analogous results for 
the minorant, 
{which were obtained by computing integrals in (\ref{eq:lower-estimate})
on each time-step cylinder.}

%From a practical point of view, the most challenging step is the reconstruction of
%$y \in Y_{\dvrg}(Q_T)$ minimizing $\maj{}$ (or $\eta \in \Ho{1}(Q_T)$ maximizing
%$\Min$).
%Since the majorant attains its minimum on $y = A \nabla u$, we reconstruct the initial
%flux using $A \nabla v$. However, since $A \nabla v$ does not necessarily belong to
%$Y_{\dvrg}(Q_T)$, we apply the smoothing operator $y = G (A \nabla u)$, where
%%
%$G:\: Y(Q_T) \rightarrow Y_{\dvrg}(Q_T)$.
% Numerically this procedure is usually inexpensive.
%--------------------------------------------------------------------------------------%
%In Table \ref{tab:example-3-maj-I-minimization-with-y}, we illustrate the minimization
%process of $\maj{}$ with respect to flux (from the left to the right) on every
%time-cylinder $Q^k$. The first and second columns contain
%the considered time layer and true error, respectively. The initial value of majorant
%$\maj{(k)} (y)$, where $y = G (A \nabla v)$, and its efficiency index are
%situated in the third and forth columns of the table. Next, we obtain $y^{opt}$ by
%minimizing the majorant, the attained minimum $\maj{(k)} (y^{opt})$
%and corresponding efficiency index form the fifth and sixth columns.
%In the final step, we minimize majorant with respect to flux constructed on the
%refined mesh, and obtain the corresponding efficiency index (see columns seven and eight).

%--------------------------------------------------------------------------------------%
\begin{table}[!t]
	\centering
	\footnotesize
	\begin{tabular}{cccccccc}
	\midrule
	$t^k$ & $\error_{(1, \: 0, \: 1)}$
			& $\maj{} \left(y\right)$             & $\Ieffmaj$ %$\Ieff \left( \maj{} \right)$
			& $\maj{} \left(y^{opt}\right)$       & $\Ieffmaj$ %$\Ieff \left( \maj{} \right)$
			& $\maj{} \left(y^{opt}_{ref}\right)$ & $\Ieffmaj$ %$\Ieff \left( \maj{} \right)$
			\\
	\midrule
	1.03 & 5.64e-04 & 4.97e-03 & 2.97 & 9.49e-04 & 1.30 & 9.22e-04 & 1.28 \\
	2.05 & 5.71e-04 & 3.87e-03 & 2.60 & 7.43e-04 & 1.14 & 7.35e-04 & 1.13 \\
	3.08 & 5.85e-04 & 3.19e-03 & 2.34 & 6.78e-04 & 1.08 & 6.75e-04 & 1.07 \\
	4.10 & 5.97e-04 & 2.73e-03 & 2.14 & 6.55e-04 & 1.05 & 6.53e-04 & 1.05 \\
	5.13 & 6.06e-04 & 2.41e-03 & 1.99 & 6.45e-04 & 1.03 & 6.45e-04 & 1.03 \\
	6.15 & 6.13e-04 & 2.17e-03 & 1.88 & 6.43e-04 & 1.02 & 6.43e-04 & 1.02 \\
	7.18 & 6.18e-04 & 1.99e-03 & 1.79 & 6.42e-04 & 1.02 & 6.42e-04 & 1.02 \\
	8.21 & 6.22e-04 & 1.84e-03 & 1.72 & 6.41e-04 & 1.02 & 6.41e-04 & 1.01 \\
	9.23 & 6.26e-04 & 1.73e-03 & 1.66 & 6.41e-04 & 1.01 & 6.41e-04 & 1.01 \\
	10.00 & 6.28e-04 & 1.65e-03 & 1.62 & 6.41e-04 & 1.01 & 6.41e-04 & 1.01 \\
	\midrule
	\end{tabular}
	\caption{Example \ref{ex:example-3}.
	Minimization of the majorant $\maj{}$ with respect to flux
	on every time-cylinder $Q^{k}$, $k = 1, ... , 39$.}
	\label{tab:example-3-maj-I-minimization-with-y}
\end{table}
\begin{table}[!t]
	\centering
	\footnotesize
	\begin{tabular}{cccccccc}
	\midrule
	$t^k$ & $\error_{(1/\sqrt{2}, \: 1/\sqrt{2}, \: 1/\sqrt{2})}$
			& $\Min^{} \left( \eta \right)$       & $\Ieffmin$ %$\Ieff$
			& $\Min^{} \left( \eta^{opt} \right)$ & $\Ieffmin$ %$\Ieff$ 
			\\
	\midrule
	1.03 & 5.64e-04 & 4.29e-04 & 0.87 & 5.09e-04 & 0.99  \\
	2.05 & 5.71e-04 & 5.11e-04 & 0.95 & 5.22e-04 & 1.00  \\
	3.08 & 5.85e-04 & 5.34e-04 & 0.96 & 5.36e-04 & 1.00  \\
	4.10 & 5.97e-04 & 5.46e-04 & 0.96 & 5.47e-04 & 1.00  \\
	5.13 & 6.06e-04 & 5.54e-04 & 0.96 & 5.55e-04 & 1.00  \\
	6.15 & 6.13e-04 & 5.60e-04 & 0.96 & 5.60e-04 & 1.00  \\
	7.18 & 6.18e-04 & 5.65e-04 & 0.96 & 5.65e-04 & 1.00  \\
	8.21 & 6.22e-04 & 5.68e-04 & 0.96 & 5.68e-04 & 1.00  \\
	9.23 & 6.26e-04 & 5.71e-04 & 0.96 & 5.71e-04 & 1.00  \\
	10.00 & 6.28e-04 & 5.73e-04 & 0.96 & 5.73e-04 & 1.00  \\
	\midrule
	\end{tabular}
	\caption{Example \ref{ex:example-3}.
	Maximization of the minorant $\Min{}$ with respect to $\eta$
	on every time-cylinder $Q^{k}$, $k = 1, ... , 39$.}
	\label{tab:example-3-min-maximization-with-eta}
\end{table}
\end{example}

%--------------------------------------------------------------------------------------%
%
%\begin{table}[!ht]
	%\centering
	%\footnotesize
	%\begin{tabular}{llllllll}
	%\midrule
	%$k$ & $[e]_{(1, \frac12)}^{(k)}$
			%& $\	MajII{(k)} (y, w)$            & $\Ieff^{(k)}$
			%& $\MajII{(k)} (y^{opt}, w^{opt})$ & $\Ieff^{(k)}$\\
	%\midrule
	%1 & 6.19e-004 & 1.81e-002 & 5.41 & 9.91e-004 & 1.27 \\
	%3 & 6.32e-004 & 6.12e-003 & 3.11 & 1.19e-003 & 1.37 \\
	%7 & 6.16e-004 & 4.48e-003 & 2.70 & 8.74e-004 & 1.19 \\
	%15 & 6.25e-004 & 2.99e-003 & 2.19 & 7.06e-004 & 1.06 \\
	%31 & 6.39e-004 & 1.93e-003 & 1.74 & 6.61e-004 & 1.02 \\
	%39 & 6.42e-004 & 1.69e-003 & 1.62 & 6.57e-004 & 1.01 \\
	%\midrule
	%\end{tabular}
	%\caption{Example \ref{ex:example-3}.
	%Minimization of majorant $\MajII{(k)}$ with respect to
	%$y$ and $w$ on every time-cylinder $Q^{k}$, where $k = 1, ... , 39$.}
	%\label{tab:example-3-maj-II-minimization-with-y-w}
%\end{table}

\begin{figure}[!h]
	\centering
	\includegraphics[trim = 0mm 0mm 0mm 1mm, clip, scale=0.9]{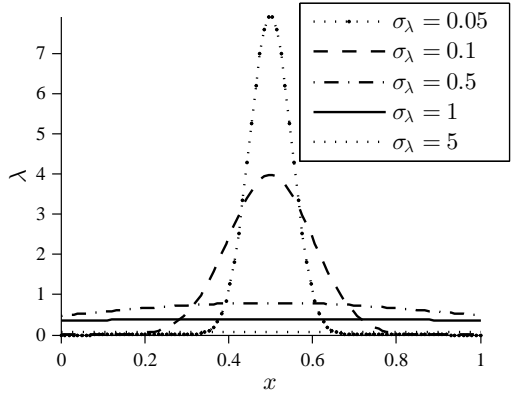}\\
	\caption{{The reaction term $\lambda(x)$ defined by (\ref{eq:gaussian-term}) }.}
	\label{fig:example-3-gaussian-function}
\end{figure}

{For the next set of numerical tests, we assume
that the reaction term is positive and behaves
as the Gaussian function, i.e., }
\begin{equation}
\lambda(x) = \frac{1}{\sigma_{\lambda} \sqrt{2 \pi}} \exp 
\left( \frac{-\left(x - \frac12 \right)^2}{2 \,\sigma_{\lambda}^2}\right).
\label{eq:gaussian-term}
\end{equation}
Then, the right-hand side of the problem is changed to 
$f = x \, (1 - x) \,(2t + 1) + (\lambda \, x \, (1 - x) + 2)\, (t^2 + t + 1)$.
From \linebreak
Fig. \ref{fig:example-3-gaussian-function}, we see that for the certain $\sigma_{\lambda}$ the reaction
$\lambda$ changes rapidly from very small values (in one part of $\Omega$) to relatively big values
(in another part).
The estimate (\ref{eq:majorant-1}) was derived specially for such type cases,
and we use this example in order to verify it.
Consider a simplified form of (\ref{eq:majorant-1}) with $\delta = 1$ and $\beta = {\rm const}$, {which implies the majorant}
\begin{multline}
	\Majmu (v, y) :=
	\| \, e (0, x) \, \|^2_{\Omega} +
	        \Int_0^T \Bigg ( \gamma \bigg\| \, \frac{\mu}{\sqrt{\lambda}} \, \R_f(v, y) \, \bigg\|^2_{\Omega} \: +
	                 \bigg(1 + \frac{1}{\beta}\bigg) \CF^2 \| \, (1 - \mu) \, \R_f(v, y) \, \|^2_{\Omega} \: +
									 (1 + \beta) \| \, \R_d(v, y) \,\|^2_{A^{ \minus 1}} \Bigg ) \dt = \\
	\| \, e (0, x) \, \|^2_{\Omega} +
	        \Int_{Q_T} \Bigg ( \gamma  \frac{\mu^2}{\lambda} \, \R^2_f(v, y) \, \: +
	                 \bigg(1 + \frac{1}{\beta}\bigg) \CF^2 (1 - \mu)^2 \, \R^2_f(v, y) \: +
									 (1 + \beta) \R^2_d(v, y) \Bigg ) \dx \dt.
	\label{eq:majorant-mu-simplified}
\end{multline}

Minimization of the right-hand side of (\ref{eq:majorant-mu-simplified}) with respect
to $\mu$ is reduced to the auxiliary variational problem: find $\hat{\mu} \in L^\infty(\Omega)$ such
that
\begin{equation}
	\Upsilon(\hat{\mu})=\!\!\inf\limits_{\mu \: \in \: L^\infty(\Omega)} \!\!\! \Upsilon(\mu), \quad \mbox{for a.e.} \quad t \in (0, T), 
	\label{eq:functional-depanding-from-mu}
\end{equation}
where
\begin{equation}
\Upsilon(\mu) := \IntQT \Bigg ( \gamma  \frac{\mu^2}{\lambda} \, \R^2_f(v, y) \, \: +
	                 \bigg(1 + \frac{1}{\beta}\bigg) \CF^2 (1 - \mu)^2 \, \R^2_f(v, y) \Bigg ) \dx \dt.
\end{equation}
It is easy to find that for a.e. $(x, t) \in Q_T$ the minimizer satisfies the condition
\begin{equation}
\hat\mu (x, t) = \frac{\CF^2(1 + \beta)\lambda}{\beta \gamma + \CF^2(1 + \beta)\lambda}.
\end{equation}
Table \ref{tab:example-3-i-eff-from-lambda} shows the efficiency of
$\Majmuhat$ and $\Majzero$ for different $\sigma_{\lambda}$. In the left part of it, the results correspond to the case where $y$ is reconstructed
by piecewise affine approximations, and in the right part we expose 
the results obtained if $y$ is taken
from a reacher space (which includes piecewise quadratic functions).
We see that $\Majzero$ grows dramatically if
$\sigma_{\lambda}$ goes to zero, while $\Majmuhat$ keeps small values of the
efficiency index for all $\sigma_{\lambda}$. In other words, $\Majmuhat$ stays efficient and robust
even if the reaction function changes its values quite drastically  in different parts of the domain.

\begin{table}[!t]
\centering
\footnotesize
\begin{tabular}{cccccc}
\midrule
$\:$      & \multicolumn{2}{c}{$y_{opt}$ (linear approximation)} & $\:$ &\multicolumn{2}{c}{$y^{ref}_{opt}$ (quadratic approximation)} \\
\cmidrule{2-3} \cmidrule{5-6}
%\midrule
$\sigma_{\lambda}$ & $\Ieffmajzero$ & $\Ieffmajmuhat$ & $\:$ &
                     $\Ieffmajzero$ & $\Ieffmajmuhat$ \\
\midrule
0.05 & $9.037 \cdot 10^8$ & 1.0080 & $\;$ &$9.037 \cdot 10^8$ & 1.0077 \\
0.10 & 2.4416             & 1.0063 & $\:$ & 2.4282            & 1.0061 \\
0.50 & 1.0020             & 1.0017 & $\:$ & 1.0020            & 1.0017 \\
1.00 & 1.0035             & 1.0029 & $\:$ & 1.0034            & 1.0027 \\
5.00 & 1.0173             & 1.0079 & $\:$ & 1.0168            & 1.0075 \\
\midrule
\end{tabular}
\caption{Example \ref{ex:example-3}. Efficiency indexes for different values of $\sigma_{\lambda}$ for $t = T$.}
\label{tab:example-3-i-eff-from-lambda}
\end{table}

%\begin{figure}[!ht]
	%\centering
	%\includegraphics[trim = 0mm 0mm 0mm 1mm, clip, scale=0.9]{pics/example-3/i-eff-mu-1-0-opt-ex-3}\\
	%\caption{Example \ref{ex:example-3}.
	%Efficiency indexes of the majorants $\Majone$, $\Majzero$, and $\Majmu$ for different $ \log \sigma_{\lambda}$.}
	%\label{fig:example-3-maj-i-eff-from-lambda}
%\end{figure}

%--------------------------------------------------------------------------------------%
\begin{example}
\rm
\label{ex:example-11}
Consider the same problem as in Example 1 but with 
$\varphi = x\, \sin\!\left(3\, \pi\, x\right)$, $A = I$, \linebreak
$\lambda(x, t) = \rho \,(t^2 + 1)(x + 10^{\minus 3})$, where $\rho$ is a positive constant,
and
%$\mathrm{e}^{t}\, \left(x \, \sin(3\, \pi \,x) \right) - ...
        %(exp(t).*(6.*pi.*cos(3.*pi.*x) - 9.*pi.^2*x.*sin(3.*pi.*x))) + ...
        %lambda(t, x, react_lb, react_ub).*(x.*sin(3.*pi.*x).*exp(t))$
$f = \mathrm{e}^{t}\, \left( x\, (1 + 9 \pi^2) \sin\!\left(3\, \pi\, x\right) -
6\, \pi\, \cos\!\left(3\, \pi\, x\right) \right) + 
\lambda (x, t) \sin( 3\, \pi \,x)\,\mathrm{e}^{t}$.
The exact solution is
$u = x\, \sin\!\left(3\, \pi\, x\right) \, \mathrm{e}^{t}$. 
%--------------------------------------------------------------------------------------%

Table \ref{tab:example-11-i-eff-from-lambda} presents the efficiency of
$\Majmuhat$, $\Majone$ and $\Majzero$ for different $\rho$.
It shows that $\Majmuhat$ always provides accurate 
upper bound of the error, whereas 
$\Majone$ and $\Majzero$ may overestimate it if $\rho$ is sufficiently small or large.
Fig. \ref{fig:example-11-i-eff-from-lambda}
illustrates the same behavior of
$\Majmuhat$, $\Majone$, and $\Majzero$ with respect to $\log \,\rho$.
These results confirm that $\Majmuhat$ is
indeed robust and 
{have serious advantages in the case where $\lambda$ may attain quite
different value (very small and very large) in different parts of the domain.}
% applicable to the cases, in which $\lambda$ changes drastically.

\begin{table}[!ht]
\centering
\footnotesize
\begin{tabular}{cccccccc}
\midrule
$\:$      & \multicolumn{3}{c}{$y_{opt}$ (linear approximation)} & $\:$ &\multicolumn{3}{c}{$y^{ref}_{opt}$ (quadratic approximation)} \\
\cmidrule{2-4} \cmidrule{6-8}
$\rho$ & $\Ieffmajone$ & $\Ieffmajzero$ & $\Ieffmajmuhat$ & $\:$ &  $\Ieffmajone$ & $\Ieffmajzero$ & $\Ieffmajmuhat$ \\
\midrule
$10^{ \minus 3}$ & 3.5751 & 94.0962 & 3.5692 & $\:$ & 2.6520 & 60.1361 & 2.6477 \\
$10^{ \minus 2}$ & 3.5835 & 29.8481 & 3.5272 & $\:$ & 2.6610 & 19.0925 & 2.6196 \\
$10^{ \minus 1}$ & 3.6209 & 9.6446 & 3.2145 & $\:$ & 2.7011 & 6.2272 & 2.4205 \\
$10^{0}$ & 3.3249 & 2.8622 & 2.1385 & $\:$ & 2.7018 & 2.1808 & 1.9152 \\
$10^{1}$ & 4.4795 & 1.2291 & 1.2608 & $\:$ & 4.2507 & 1.2200 & 1.2608 \\
$10^{2}$ & 11.7222 & 1.0795 & 1.0787 & $\:$ & 11.2044 & 1.0795 & 1.0787 \\
$10^{3}$ & 37.6760 & 1.0622 & 1.0622 & $\:$ & 36.4553 & 1.0622 & 1.0622 \\
$10^{3}$ & 37.6760 & 1.0622 & 1.0622 & $\:$ & 36.4553 & 1.0622 & 1.0622 \\
\midrule
\end{tabular}
\caption{Example \ref{ex:example-11}. Efficiency indexes for different values of $\rho$ 
{for $t = T$}.}
\label{tab:example-11-i-eff-from-lambda}
\end{table}

\begin{figure}[!ht]
	\centering
	\includegraphics[trim = 0mm 0mm 0mm 1mm, clip, scale=0.9]{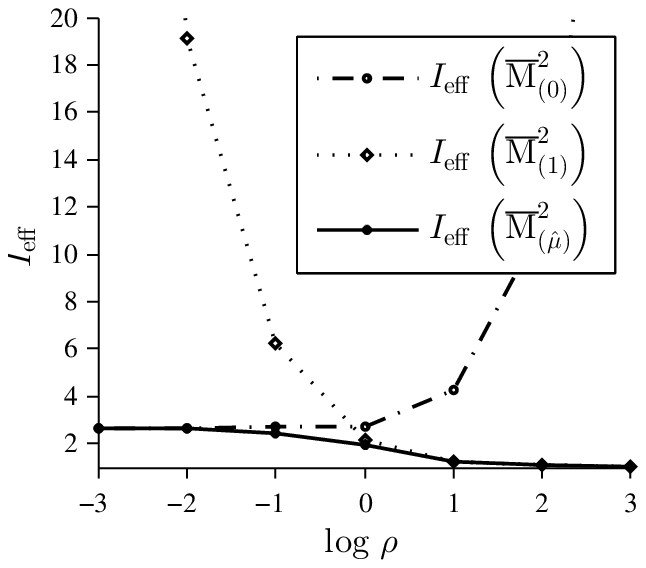}\\
	\caption{Example \ref{ex:example-11}.
	Efficiency indexes of the majorants $\Majone$, $\Majzero$, and $\Majmuhat$ for different constant $ \log \rho$.}
	\label{fig:example-11-i-eff-from-lambda}
\end{figure}

\end{example}
%
%--------------------------------------------------------------------------------------%
% example-33
%--------------------------------------------------------------------------------------%

\begin{example}
\rm
\label{ex:example-33}
%
%--------------------------------------------------------------------------------------%
%
We set $\Omega = (0, 1)$, $T$ = 10, and consider the problem with mixed 
Dirichlet--Neumann boundary conditions, namely, $u(0, t) = 0$ and 
$u^{\prime}(1, t) = 0$, 
and 
\begin{align}
	f & \; = x\, \Big(\cos\!\left(t\right) - \sin\!\left(t\right)\Big)\, {\left(x - 1\right)}^2 -
	             2 \Big(\cos\!\left(t\right) + \sin\!\left(t\right) \Big) \, \Big( 3\, x - 2 \Big), \nonumber \\
%	g & \; = \Big(t\, \cos\!\left(t\right) + 1 \Big)\, \bigg( \Big(\cos\!\left(\pi\, x\right) + 1\Big)\, \Big( \sin\!\left(\pi\, x\right) +  \pi\, \cos\!\left(\pi\, x\right) \Big)	
%	         - \pi\, {\sin\!\left(\pi\, x\right)}^2 \bigg), \nonumber \\
\varphi & \; = \sin\!\left(\pi\, x\right)\, \Big(\cos\!\left(\pi\, x\right) + 1\Big),
\nonumber\\ 
{ \lambda} & \; = 0.
	\label{eq:ex-33}
\end{align}
The exact solution is
$u(x, t) = \sin\!\left(\pi\, x\right)\, \left(t\, \cos\!\left(t\right) + 1\right)\,
\left(\cos\!\left(\pi\, x\right) + 1\right)$.
%
%\begin{figure}[!ht]
	%\centering
	%\includegraphics[trim = 0mm 2mm 0mm 5mm, clip=true, scale=1]{pics/example-33/ex-33-approximation-v}
	%\caption{{Example \ref{ex:example-33}}.
	%The approximate solution $v$ constructed on the mesh $\Theta_{40 \times 40}$.}
	%\label{fig:example-33-approximation-v}
%\end{figure}
%
%--------------------------------------------------------------------------------------%

Table \ref{tab:example-33-error-maj-ieff-for-different-meshes} shows the
relative error $\frac{\error_{(1, \: 0, \: 1)}}{[u]^2}$, the majorant
$\frac{\maj{}}{[u]^2}$, and the efficiency index for different meshes, where $N_1$
denotes the amount of intervals with respect to the space coordinate $x$ and $K$ with respect to the time $t$.
We can see that the efficiency index stays on the approximately
same level for all considered meshes, therefore
the majorant does not deteriorate in the process 
of mesh refining. 
It is worth remarking that results exposed in Table 
\ref{tab:example-33-error-maj-ieff-for-different-meshes} are quite typical, and similar behavior of the error majorant was observed in many other numerical tests.
\begin{table}[!ht]
	\centering
	\footnotesize
	\begin{tabular}{ccccc}
	\midrule
	$N_1$ & $K$ & $\error_{(1, \: 0, \: 1)} \,/\, [u]^2$ &
	$\maj{} \,/\, [u]^2$ & $\Ieffmaj$\\
	\midrule
	20 & 20  & 6.41e-03 & 2.23e-02 & 1.87 \\
	20 & 40  & 5.87e-03 & 2.27e-02 & 1.96 \\
	20 & 80  & 5.88e-03 & 2.26e-02 & 1.96 \\
	20 & 160 & 5.89e-03 & 2.27e-02 & 1.96 \\
	\midrule
	40 & 40  & 1.37e-03 & 5.14e-03 & 1.94 \\
	40 & 80  & 1.34e-03 & 5.04e-03 & 1.94 \\
	40 & 160 & 1.34e-03 & 5.14e-03 & 1.95 \\
	\midrule
	80 & 80  & 3.26e-04 & 1.22e-03 & 1.93 \\
	\midrule
	\end{tabular}
	\caption{{Example \ref{ex:example-33}}.
	The relative error, majorant, and its efficiency index with respect to different 	
	meshes $\Theta{\:K \times N_1}$ for $t = T$.}
	\label{tab:example-33-error-maj-ieff-for-different-meshes}
\end{table}

%--------------------------------------------------------------------------------------%
\begin{figure}[!t]
	\centering
	\includegraphics[scale=0.9]{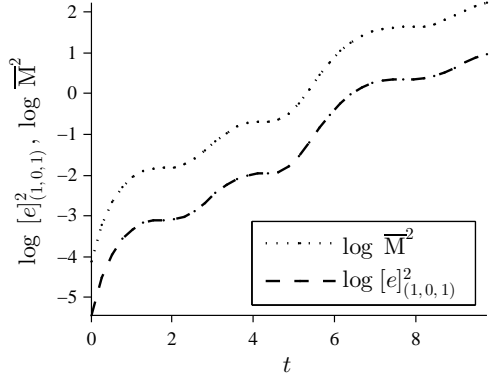}
	\caption{{Example \ref{ex:example-33}}.
	The logarithm of the error and majorant with respect to time.}
	\label{fig:example-33-log-maj-error}
\end{figure}

In Fig. \ref{fig:example-33-log-maj-error}, we depict growth of the error and majorant in the logarithmic scale (for mesh 
$\Theta_{\:K \times N_1}$ = $\Theta_{\:40 \times 40}$). The gap between two curves 
reflects the impact of the term ${\overline{\mathrm m}}^2_{\mathrm{f}}$ (see Table 
\ref{tab:example-33-majorant-terms}) and corresponds to the efficiency index 
$\Ieffmaj  \approx 1.87, \ldots, 1.96$. The efficiency of the majorant can be 
improved by choosing higher order approximations for the flux reconstruction
%to reconstruct the flux $y$,  adding 
%element-wise based bubble functions 
(e.g., we can add element-wise based bubble functions; corresponding results are presented in columns 7 and 8 of Table \ref{tab:example-33-maj-minimization-with-y}).
%--------------------------------------------------------------------------------------%
Table \ref{tab:example-33-maj-minimization-with-y} illustrates minimization of $\maj{}$ with respect to 
$y$ and corresponding efficiency indexes on every local time-cylinder $Q^k$. 
For example, consider the row of Table \ref{tab:example-33-maj-minimization-with-y}
related to $t_k = 6.15$. The corresponding error $\error_{(1, \: 0, \: 1)} = 1.33e
\minus 03$
and the estimate provided by the majorant with a simple (patch averaging) reconstruction
 of the flux is $\maj = 6.42e\minus 03$ $\left( \Ieffmaj = 2.20 \right)$. If we apply a more sophisticated procedure
and reconstruct flux by minimizing the majorant with respect to values of $y$ associated 
with the given spatial mesh on the time-layer $t_k = 6.15$, then we obtain the efficiency
index $\Ieffmaj = 2.18$. If we use a twice finer spacial mesh, then the index 
decreases up to $1.90$.  Obtained results show that for this particular example simple flux reconstructions generate sufficiently accurate estimates. 

It is worth noting that, in general, 
the efficiency of the majorant depends on the several factors. First, it can be improved
by using advanced reconstructions of the flux (e.g., by adding extra degrees of 
freedom, bubble functions, etc). However, this approach may lead to a limited effect,
if approximations are coarse with respect to the time variable and/or the time steps 
are large (the term $v_t$ in the balancing term may be a piecewise 
constant function, which cannot reflect the behavior of $f$).
% 
%Consider $Q^k$, which corresponds to the time $t_k = 6.15$. The minimization of the majorant with respect to 
%$y$ (represented by linear approximations in time) provides the improvement of the  efficiency index from $\Ieffmaj = 2.2$ to $\Ieffmaj = 2.18$. However, if we add the
%bubble function into the approximation of $y$ and run the optimization process again,
%we obtained the improved value of the efficient index $\Ieff = 1.90$. Therefore, 
%increasing the order of approximation functions for the flux indeed improves the 
%accuracy of the error estimate.
%
\begin{table}[!ht]
	\centering
	\footnotesize
	\begin{tabular}{ccccc}
	\midrule
	$t^k$ & ${{\mid\mid\mid \! \nabla(u - v) \! \mid\mid\mid}^2_{ \overline{Q}_{{t \,}^k} }} $
			& $\| \: (u - v)(x, t^k) \: \|^2_{\Omega}$
			& ${{\| \: y - \nabla v \: \|}^2_{ \overline{Q}_{{t \,}^k} }}$
			& ${{\| \: f - v_t + \dvrg y \: \|}^2_{ \overline{Q}_{{t \,}^k} }}$ \\
	\midrule
	1.03 & 8.72e-03 & 1.42e-04 & 8.37e-03 & 7.48e-02 \\
	2.05 & 6.00e-04 & 4.40e-06 & 5.02e-04 & 7.71e-03 \\
	3.08 & 1.29e-02 & 2.15e-04 & 1.24e-02 & 1.11e-01 \\
	4.10 & 1.12e-02 & 1.62e-04 & 1.01e-02 & 1.02e-01 \\
	5.13 & 2.07e-02 & 3.16e-04 & 2.07e-02 & 1.87e-01 \\
	6.15 & 1.62e-01 & 2.53e-03 & 1.60e-01 & 1.38e+00 \\
	7.18 & 1.30e-01 & 2.00e-03 & 1.27e-01 & 1.11e+00 \\
	8.21 & 3.36e-03 & 3.14e-05 & 3.35e-03 & 5.98e-02 \\
	9.23 & 2.05e-01 & 3.10e-03 & 2.00e-01 & 1.75e+00 \\
	10.00 & 2.22e-01 & 3.33e-03 & 2.14e-01 & 1.90e+00 \\
	\midrule
	\end{tabular}
	\caption{{Example \ref{ex:example-33}}.
	Two terms of the error and two terms of the majorant with respect to time.}
	\label{tab:example-33-majorant-terms}
\end{table}
%
%--------------------------------------------------------------------------------------%

\begin{table}[!ht]
	\centering
	\footnotesize
	\begin{tabular}{cccccccc}
	\midrule
	$t^k$ & $\error_{(1, \, 0, \, 1)}$
			& $\maj{}(y)$             & $\Ieffmaj$
			& $\maj{}(y^{opt})$       & $\Ieffmaj$
			& $\maj{}(y^{opt}_{ref})$ & $\Ieffmaj$ \\
	\midrule
	1.03 & 1.33e-03 & 6.58e-03 & 2.22 & 6.50e-03 & 2.21 & 4.79e-03 & 1.90 \\
	2.05 & 1.46e-03 & 5.31e-03 & 1.91 & 5.31e-03 & 1.91 & 5.25e-03 & 1.90 \\
	3.08 & 1.38e-03 & 6.10e-03 & 2.11 & 6.05e-03 & 2.10 & 4.98e-03 & 1.90 \\
	4.10 & 1.47e-03 & 5.74e-03 & 1.98 & 5.73e-03 & 1.97 & 5.27e-03 & 1.89 \\
	5.13 & 1.42e-03 & 5.93e-03 & 2.04 & 5.89e-03 & 2.04 & 5.15e-03 & 1.91 \\
	6.15 & 1.33e-03 & 6.42e-03 & 2.20 & 6.36e-03 & 2.18 & 4.81e-03 & 1.90 \\
	7.18 & 1.40e-03 & 5.65e-03 & 2.01 & 5.63e-03 & 2.01 & 5.03e-03 & 1.90 \\
	8.21 & 1.43e-03 & 5.16e-03 & 1.90 & 5.16e-03 & 1.90 & 5.15e-03 & 1.90 \\
	9.23 & 1.39e-03 & 5.67e-03 & 2.02 & 5.65e-03 & 2.02 & 5.01e-03 & 1.90 \\
	10.00 & 1.41e-03 & 5.58e-03 & 1.99 & 5.56e-03 & 1.99 & 5.06e-03 & 1.90 \\
	\midrule
	\end{tabular}
	\caption{Example \ref{ex:example-33}.
	Minimization of the majorant $\maj{}$ with respect to
	flux on every time-cylinder $Q^{k}$, $k = 1, ... , 39$.}
	\label{tab:example-33-maj-minimization-with-y}
\end{table}
%
%--------------------------------------------------------------------------------------%

%Another important question, which should be discussed in the context of a posteriori
%error estimation concerns indication of the local error.  Below we show that the
%majorant implies such an error indicator that can be used for restructuring of the
%space-time mesh in the process of integration. First, we note that the majorant can be
%represented in the form
%%
%\begin{equation}
	%\maj{} (v, y, \beta) : =
	%\| \, e(x, 0) \, \|^2_{\Omega} +
	%\frac1\delta \Int_0^T
	%\Bigg(
	%\big(1 + \beta \big) \mdI +
	%\bigg(1 + \frac{1}{\beta}\bigg) \CF^2 {\overline{\mathrm m}}^2_{\mathrm{f}}
	%\Bigg) \dt,
	%\label{eq:indicator-components}
%\end{equation}	
%%
%where
%%
%\begin{equation}
	%\mdI{} = 	\| \: y - \nabla v \: \|^2_{\Omega}, \quad
	%{\overline{\mathrm m}}^2_{\mathrm{f}} = \| \: f - v_t + \dvrg y \: \|^2_{\Omega}.
	%\label{eq:indicator-components-2}
%\end{equation}
%%
%Here, ${\overline{\mathrm m}}^2_{\mathrm{f}}$ is `reliability term', which is necessary
%to provide a guaranteed upper bound, and the major part of the error is usually encompassed in $\mdI{}$.
%It is natural to use the latter term as an error indicator. 
%By numerical results presented in
%Fig. \ref{fig:example-33-error-maj-time}, we confirm that $\mdI{}$ generates the efficient indicator for the
%error $\error_{(\sqrt{2-\delta}, \, 0, \, 1)}$ over $\Omega$ in every time-layer.
%
%
Now, we shortly discuss results related to error indicators generated by the majorant,
which can be represented in the form
\begin{equation}
	\maj{} (v, y, \beta) : =
	\| \, e(x, 0) \, \|^2_{\Omega} +
	\frac1\delta \Int_0^T
	\Bigg(
	\big(1 + \beta \big) \mdI +
	\bigg(1 + \frac{1}{\beta}\bigg) \, \CF^2 \mfI
	\Bigg) \dt,
	\label{eq:indicator-components}
\end{equation}	
where
\begin{equation}
	\mdI = \| \: y - \nabla v \: \|^2_{\Omega}, \quad
	\mfI = \| \: f - v_t - \lambda v + \dvrg y \: \|^2_{\Omega}.
	\label{eq:indicator-components-2}
\end{equation}
The `reliability term' $\mfI$ is necessary to provide a guaranteed upper 
bound, but the major part of the error is usually encompassed in $\mdI$. 
%Therefore, it 
%is natural to use the latter term as an error indicator. 
%
The conclusion, which follows from our experience is that the term $\mdI$ (which is an
integral over $\Omega$) can be considered as an efficient indicator of element-wise 
error. 
Fig. \ref{fig:example-33-error-maj-time} presents distribution of the error indicator
$\mdI$ and error over time-step cylinder $Q_{30}$ (for zero reaction function).
Assume now that \linebreak 
$\lambda(x) = \frac{1}{\sigma_{\lambda} \sqrt{2 \pi}} 
\exp \left(\frac{-\left(x - \frac12 \right)^2}{2 \,\sigma_{\lambda}^2}\right)$
(which correspondingly changes $u$ and $f$) 
with
$\sigma_{\lambda}$ = 0.05 and $\sigma_{\lambda}$ = 0.1.
Fig. \ref{fig:example-33-error-maj-distr-with-reaction-from-time} illustrates typical results related to different $\sigma_{\lambda}$	and time-cylinder $Q_{30}$. 

%Numerical results depicted in Fig. \ref{fig:example-33-error-maj-time} (and similar results obtained for other 
%tests) allow us to conclude that $\mdI{}$ generates the efficient error indicator. The problem was also tested for
% From Fig. \ref{fig:example-33-error-maj-distr-with-reaction-from-time}, 
%we see that efficiency of $\mdI{}$ does not deteriorate even with drastically changing reaction $\lambda$.

\begin{figure}[!t]
	\centering
	% ${m_d}^{\oplus}_{\mathrm{I}}$
	%\subfloat[k = 30]{
	\includegraphics[scale=0.9]{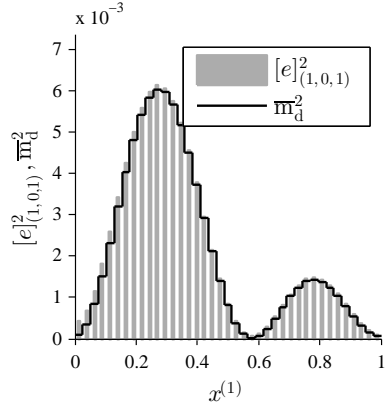}
	%}
	%\subfloat[k = 40]{
	%\includegraphics[scale=0.9]{ex-33-maj-I-opt-y-distr-e-maj-time-40}}
	\caption{{Example \ref{ex:example-33}}.
	The distribution of the local errors and indicator $\mdI{}$ for time-layers $Q^{30}$.}
	\label{fig:example-33-error-maj-time}
\end{figure}

\begin{figure}[!t]
	\centering
	% ${m_d}^{\oplus}_{\mathrm{I}}$
	\subfloat[$\sigma_{\lambda}$ = 0.05, $Q^{30}$]{
	\includegraphics[scale=0.9]{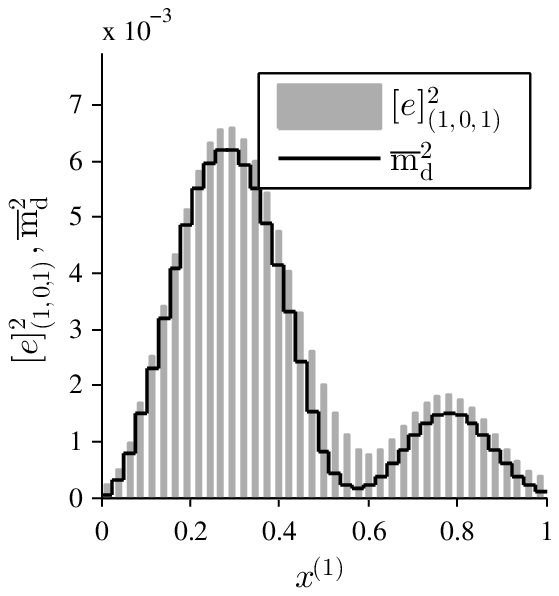}}
	%\subfloat[$\sigma_{\lambda}$ = 0.05, k = 40]{
	%\includegraphics[scale=0.9]{ex-33-maj-I-mu-opt-sigma-5-1e-2-distr-e-maj-time-40}}\\
	\subfloat[$\sigma_{\lambda}$ = 0.1, $Q^{30}$]{
	\includegraphics[scale=0.9]{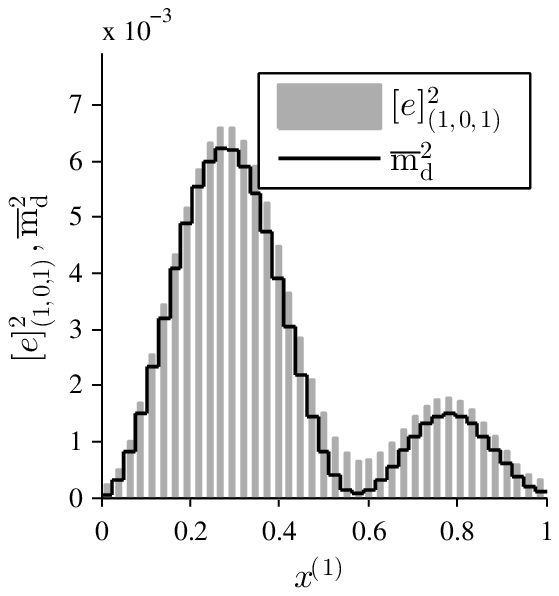}}
	%\subfloat[$\sigma_{\lambda}$ = 0.1, k = 40]{
	%\includegraphics[scale=0.9]{ex-33-maj-I-mu-opt-sigma-1e-1-distr-e-maj-time-40}}\\
	%\subfloat[$\sigma_{\lambda}$ = 0.5, k = 30]{
	%\includegraphics[scale=0.9]{ex-33-maj-I-mu-opt-sigma-5-1e-1-distr-e-maj-time-30}}
	%\subfloat[$\sigma_{\lambda}$ = 0.5, k = 40]{
	%\includegraphics[scale=0.9]{ex-33-maj-I-mu-opt-sigma-5-1e-1-distr-e-maj-time-40}}
	\caption{{Example \ref{ex:example-33}}.
	The error and indicator for time-layers $Q^{30}$ for $\sigma_{\lambda}$ = 0.05 and $\sigma_{\lambda}$ = 0.1.}
	\label{fig:example-33-error-maj-distr-with-reaction-from-time}
\end{figure}

\end{example}

%--------------------------------------------------------------------------------------%
% example-8-2d
%--------------------------------------------------------------------------------------%

\begin{example}
\label{ex:example-8-2d}
\rm
Consider the problem $\Omega = [0, 1] \times [0, 1]\in {\mathds{R}}^2$,  
$T = 1$, and homogeneous Dirichlet boundary condition, 
$\varphi(x) = \sin\!\left(\pi\, x\right)\, \sin\!\left(3\, \pi\, y\right) +
\sin\!\left(3\, \pi\, x\right)\, \sin\!\left(\pi\, y\right)$, $A = I$,  and
\begin{equation}
	f =  \Big(\sin\!\left(\pi\, x\right)\, \sin\!\left(3\, \pi\, y\right) +
	\sin\!\left(3\, \pi\, x\right)\, \sin\!\left(\pi\, y\right)\Big)\,
	\Big(\cos\!\left(t\right) + 10\, {\pi}^2\, \sin\!\left(t\right) +
	10\, {\pi}^2\, t^3 + 3\, t^2 + 10\, {\pi}^2 \Big).
	\label{eq:ex-8-right-hand-side}
\end{equation}
The corresponding exact solution
$u = \Big(\sin(\pi x) \sin(3 \pi y) + \sin(\pi y) \sin(3 \pi x) \Big)
	    \Big(t^3 + \sin(t) + 1 \Big)$
is a rapidly changing function.
%
%\begin{equation}
%	u = \Big(\sin(\pi x) \sin(3 \pi y) + \sin(\pi y) \sin(3 \pi x) \Big)
%	    \Big(t^3 + \sin(t) + 1 \Big).
%	\label{eq:ex-8-exact-solution}
%\end{equation}
%
%--------------------------------------------------------------------------------------%
%

In Table \ref{tab:error-maj-i-eff-ex-8}, we compare 
$\frac{\error_{(1, \, 0, \, 1)}}{[u]^2}$ with
$\frac{\maj{}}{[u]^2}$ 
and its efficiency index (for the approximate solution computed
on the mesh $\Theta_{\:K \times N_1 \times N_2}$ = 
$\Theta_{\:100 \times 50 \times 50}$). We can see that the majorant based on 
{`cheap' (local patch averaging) reconstruction of the flux 
$y = G(\nabla v)$ 
provides a quite realistic upper bound of the error. 
However, 
in more complicated problems an optimization of the majorant with respect to $y$ 
may be useful. This procedure
yields sharper upper 
bounds but requires more computational efforts (concerning the corresponding methods
based on multigrid, isogeometric elements, and other methods, see 
\cite{Valdman2009, KleissTomar2012, MaliNeittaanmakiRepin2013}). 
%The optimization of the majorant with respect to $y$ requires usually more time 
%resources, therefore highly efficient algorithm mush be applied (the research on 
%this topic can be found in \cite{Valdman2009, KleissTomar2012}). 
}	
{
\begin{table}[!t]
	\centering
	\footnotesize
	\begin{tabular}{cccc}
	\midrule
	$t^k$ & $\error_{(1, \, 0, \, 1)} / [u]^2$
			& $\maj{} (G(\nabla v)) / [u]^2 $ & $\Ieffmaj$ \\
			%& $\Ieff \left( \maj{} \right)$
			%& $\maj{} (\nabla u)$ & $\Ieff \left( \maj{} \right)$\\
	\midrule
	0.10 & 2.58e-03 & 2.56e-02 & 3.15\\
	0.20 & 2.70e-03 & 2.68e-02 & 3.15\\
	0.30 & 2.74e-03 & 2.72e-02 & 3.15\\
	0.40 & 2.76e-03 & 2.74e-02 & 3.15\\
	0.51 & 2.78e-03 & 2.75e-02 & 3.15\\
	0.60 & 2.79e-03 & 2.76e-02 & 3.15\\
	0.70 & 2.79e-03 & 2.77e-02 & 3.15\\
	0.80 & 2.80e-03 & 2.78e-02 & 3.15\\
	0.90 & 2.81e-03 & 2.79e-02 & 3.15\\
	1.00 & 2.82e-03 & 2.80e-02 & 3.15\\
	\midrule
	\end{tabular}
	\caption{Example \ref{ex:example-8-2d}.
	The relative true error and relative majorant with respect to time.}
	\label{tab:error-maj-i-eff-ex-8}
\end{table}
}
Next goal is to investigate the accuracy of the error indicator defined in
(\ref{eq:indicator-components})--(\ref{eq:indicator-components-2}). We analyze two 
different measures, which can be called `weak' and `strong' and are discussed in details in \cite{MaliNeittaanmakiRepin2013}. The first measure is 
studied in the context of a certain marking procedure $\Marker$, 
which maps element-wise error into a boolean array, i.e., it deals with 
the values $0$ and $1$ only. 
The corresponding `weak' measure $\mathcal{M}_{\rm weak} \in [0, 1]$ is defined by the 
percentage of correctly marked elements.
%For example, if marking based on true error coincides with the marking
%based on error indicator, then

%compares
%markers by the truth error and marker by the indicator generated by the majorant.
%Therefore, the indicator is $\eps_{\rm weak}$-accurate with
%respect to the marker $\Marker$ if the following inequality holds
%%
%\begin{equation}
	%\mathcal{M}_{\rm weak} (\mdI{}) :=
	%1-\frac{\bnorm{ \Marker (\error_{(1, \: 0, \: 1)}) \equiv \Marker (\mdI{}) }}{N} \;
	%\leq\; \eps_{\rm weak},
	%\label{eq:weak-accuracy-definition}
%\end{equation}
%%
%where $\bnorm{ a}$ stands for the sum $\Sum^N_{i=1} a_i$, and $\equiv$ denotes the
%logical equivalence rule.

Another measure compares element-wise values of the true error and estimates of local
errors generated by the error indicator.
For $\mdI{}$, it is defined by the relation
\begin{equation}
	\mathcal{M_{\rm str}} (\mdI{}) :=
	\frac{ \big| \error_{(1, \: 0, \: 1)} - \mdI{} \big|}{ \big| \error_{(1, \: 0, \: 1)} \big|}.
\end{equation}

To analyze the  quality of the weak measure, we consider bulk marking procedure 
$\Marker_{ \theta}$, where $\theta \in (0, 1)$ (see, e.g., \cite{Dorfler1996}).
Fig.
\ref{fig:error-maj-bulk-marking-ex-8} illustrates $\Marker_{ \theta}$ for 
$\theta$ = 0.2 and 0.4 (it has been performed for the 
actual error and for the error indicator $\mdI{}$). 
The results obtained for the error 
(Fig. \ref{fig:error-maj-bulk-marking-ex-8} left) and for the indicator (right) are almost identical. 
We also have obtained quite small values of the weak error measure for bulk 
parameters $\theta$ = 0.2, 0.3, and 0.4:
\begin{equation}
\mathcal{M_{\rm weak}}(\Marker_{ 0.2}) = 5.83e \minus 03, \quad
\mathcal{M_{\rm weak}}(\Marker_{ 0.3}) = 4.58e \minus 03, \quad \mbox{and} \quad
\mathcal{M_{\rm weak}}(\Marker_{ 0.4}) = 1.41e \minus 02.
\end{equation}
\begin{figure}[!t]
	\centering
	\subfloat{\includegraphics[trim = 0mm 0mm 0mm 1mm, clip, scale=1]{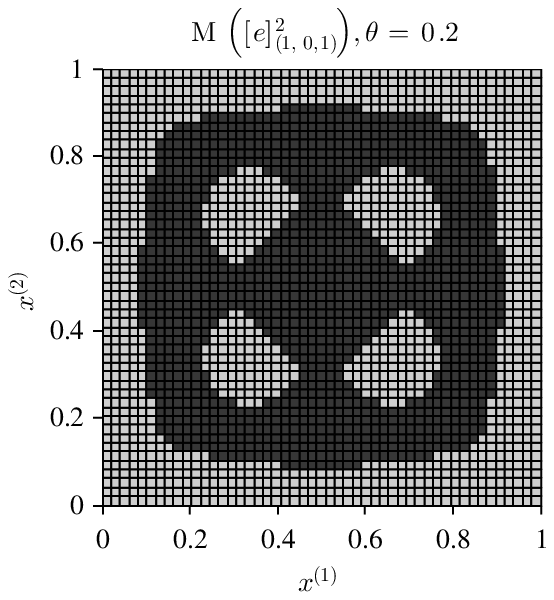}}
	\subfloat{\includegraphics[trim = 0mm 0mm 0mm 1mm, clip, scale=1]{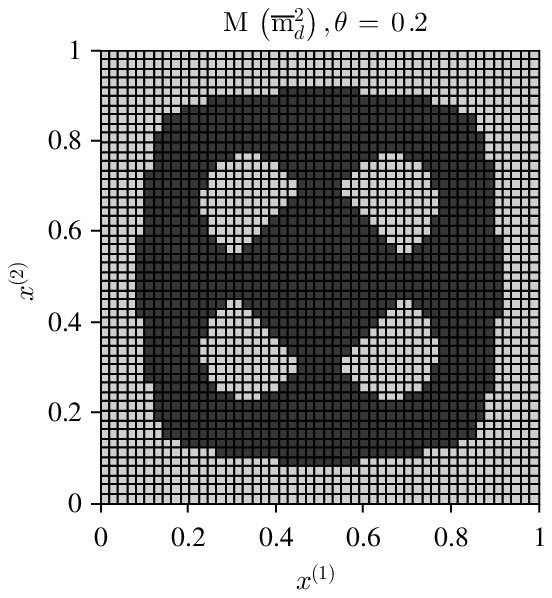}}\\
	\vskip -6pt
	\subfloat{\includegraphics[trim = 0mm 0mm 0mm 1mm, clip, scale=1]{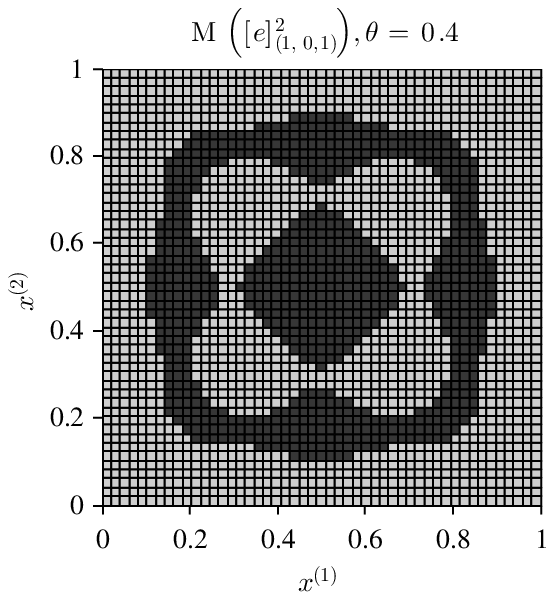}}
	\subfloat{\includegraphics[trim = 0mm 0mm 0mm 1mm, clip, scale=1]{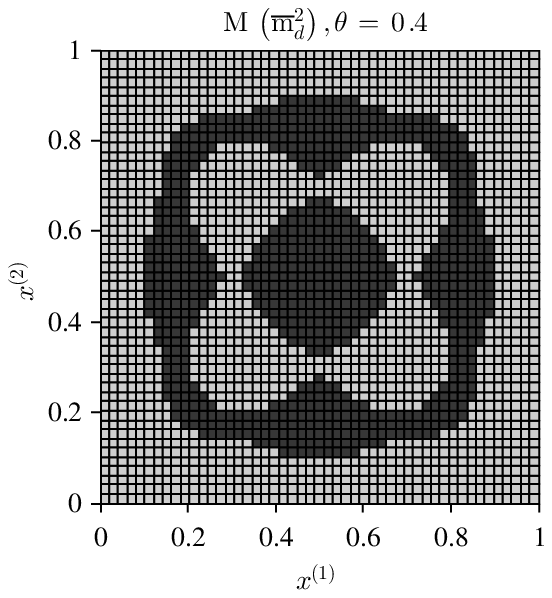}}\\
	\caption{Example \ref{ex:example-8-2d}.
	'Bulk' marking for $\theta$ = 0.2 and 0.4 based on the true error $\error_{(1, \, 0, \, 1)}$ (left)
	and the indicator $\mdI{} (G(\nabla v))$ for $Q^{40}$.}
	\label{fig:error-maj-bulk-marking-ex-8}
\end{figure}
To understand whether or not the error indicator is quantitatively sharp and reproduces
the error distribution accurately, we consider the histograms depicted in Fig. 
\ref{fig:error-maj-hist-ex-8}, which are constructed by the procedure suggested in 
\cite{MaliNeittaanmakiRepin2013}. We assume that all element-wise errors are ranked in the decreasing order with respect to values of the true error distribution $\error_{(1, \, 0, \, 1)}$, and renumber all the elements accordingly, so that the element with the
 largest error is numbered 1. Then, we depict errors in this new order. 
The distribution of the element-wise errors generated by $\mdI{}$ is depicted in the same way.
In Fig. \ref{fig:error-maj-hist-ex-8}, we consider true and
indicated error distributions (the approximation computed on a regular mesh with 
2 500 elements). If $\mdI{}$ is accurate in the strong sense, then the corresponding histogram (on the right) must resemble the histogram generated by the true error.
Therefore, Fig. \ref{fig:error-maj-hist-ex-8} shows that in this example the indicator is 
indeed sharp in the strong sense.

%It is depicted starting from the element with the largest error and ending with the element with
%the smallest error. The order of elements exposed along the horizontal axis is fixed and error
%indicator distribution $\mdI{}$ is presented in the same order. Therefore, 

\begin{figure}[!t]
	\centering
	\vskip -5pt
	\subfloat{\includegraphics[trim = 0mm 0mm 0mm 1mm, clip, scale=1]{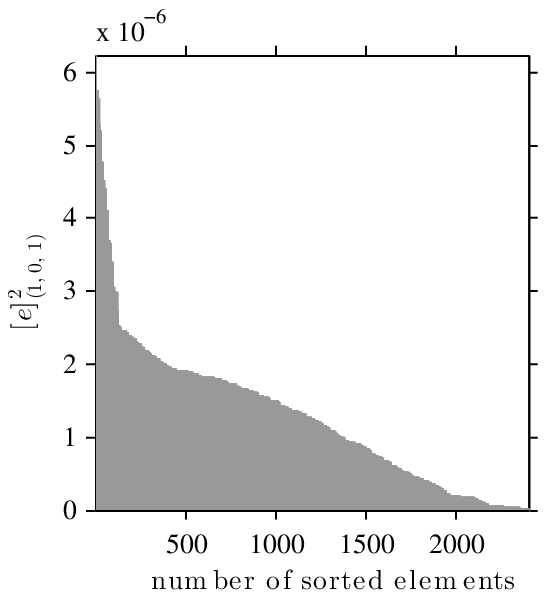}}
	\subfloat{\includegraphics[trim = 0mm 0mm 0mm 1mm, clip, scale=1]{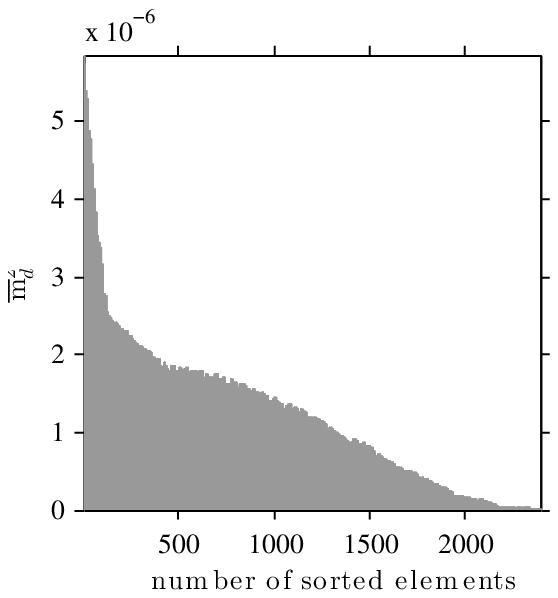}}\\
	\caption{Example \ref{ex:example-8-2d}.
	Histograms of the ranked element-wise errors and indicator by $\mdI{}$
	for $Q^{40}$ (2 500 elements).}
	\label{fig:error-maj-hist-ex-8}
\end{figure}

\end{example}

%--------------------------------------------------------------------------------------%
% example-14-2d
%--------------------------------------------------------------------------------------%

\begin{example}
\label{ex:example-14-2d}
\rm
%--------------------------------------------------------------------------------------%
%
Finally, we test the example, in which the exact solution essentially changes both in space and time.
%(see Fig. \ref{fig:approx-sol-ex-14}). 
Here,  
$\Omega = [0, 1] \times [0, 1]$, $T = 10$, $\partial \Omega = \Gamma_D$, $u = 0$ on
$S_D$, and  
$\varphi(x) = \sin\!\left(\pi\, x\right)\, \sin\!\left(2\, \pi\, y\right) +
\sin\!\left(2\, \pi\, x\right)\, \sin\!\left(\pi\, y\right)$, 
$A = I$, and
\begin{multline}
	f = \sin\!\left(2\, \pi\, x\right)\, \sin\!\left(\pi\, y\right)\,
	    \bigg(\cos\!\left(t\right) - \sin\!\left(t\right) + 
			5\, {\pi}^2\, \Big( \cos\!\left(t\right) + \sin\!\left(t\right) \Big) \bigg) \, + \\
	  	\sin\!\left(\pi\, x\right)\, \sin\!\left(2\, \pi\, y\right)\,
	    \bigg(\sin\!\left(t\right) + t\, \cos\!\left(t\right) + 
			5\, {\pi}^2\, \Big(t\, \sin\!\left(t\right) + 1 \Big) \bigg).
	\label{eq:ex-14-2d-right-hand-side}
\end{multline}
%
%--------------------------------------------------------------------------------------%
The exact solution is $u = \sin\!\left(\pi\, x\right)\, \sin\!\left(2\, \pi\, y\right)\,
\left(t\, \sin\!\left(t\right) + 1\right) + \sin\!\left(2\, \pi\, x\right)\,
\sin\!\left(\pi\, y\right)\, \left(\cos\!\left(t\right) + \sin\!\left(t\right)\right)$.

We investigate the efficiency of bulk marking for $\theta$ = 0.5 and different 
time-layers for approximate solution constructed on the mesh 
$\Theta_{\:K \times N_1 \times N_2} = \Theta_{\: 200 \times 50 \times 50}$ (see Fig.
\ref{fig:error-maj-bulk-marking-y-ux-ex-14}). Here, 
\begin{equation*}
\mathcal{M_{\rm weak}}(\Marker_{ \rm 0.5}) = 5.08e \minus 02\;\; \mbox{for}\;\; Q^{40}, \quad
\mathcal{M_{\rm weak}} (\Marker_{ \rm 0.5}) = 5.4e \minus 02\;\; \mbox{for}\;\; Q^{80}, \quad
\mbox{and} \;\;
\mathcal{M_{\rm weak}}(\Marker_{ \rm 0.5}) = 3.70e \minus 02\;\; \mbox{for}\;\; Q^{160}.
\end{equation*}
The histograms from Fig. 
\ref{fig:error-maj-hist-ex-14} are
constructed by the same method as in the previous example (for 2 500 elements).
They confirm quantitative efficiency of $\mdI{}$ .  
Again, we see that the majorant provides an efficient estimation of the global error 
as well as indication of element-wise errors. In this and also
other experiments, we have observed that the quality of error estimation is good
if the solution is smooth in space and rather monotonic in time, and it becomes less
accurate if the solution admits large gradients with respect to the spacial variables
and large time derivatives.
In our experiments, this was mainly due to rather coarse (piecewise affine) 
approximations
in time used for $v$ and $y$. This fact constrains accuracy of the error estimation 
(especially in the context of the reliability term $\mfI$).
However, choosing richer spaces for the reconstruction of $v$ and $y$ will lead to sharper estimates. 

\begin{figure}[!ht]
	\centering
	\subfloat[$t^k$ = 2.06]{\includegraphics[trim = 0mm 0mm 0mm 2mm, clip=true, scale=1]{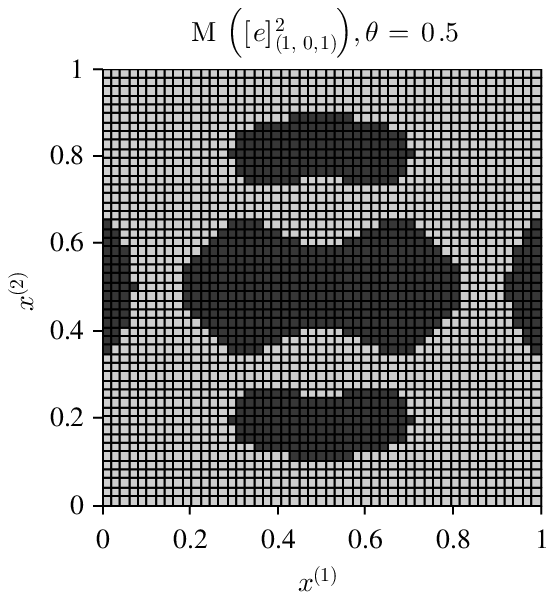}}
	\subfloat[$t^k$ = 2.06]{\includegraphics[trim = 0mm 0mm 0mm 2mm, clip=true, scale=1]{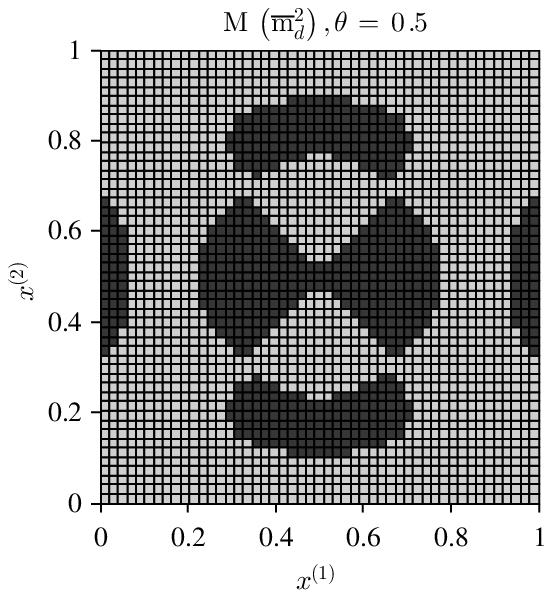}}\\
	\vskip -4pt
	\subfloat[$t^k$ = 4.07]{\includegraphics[trim = 0mm 0mm 0mm 2mm, clip=true, scale=1]{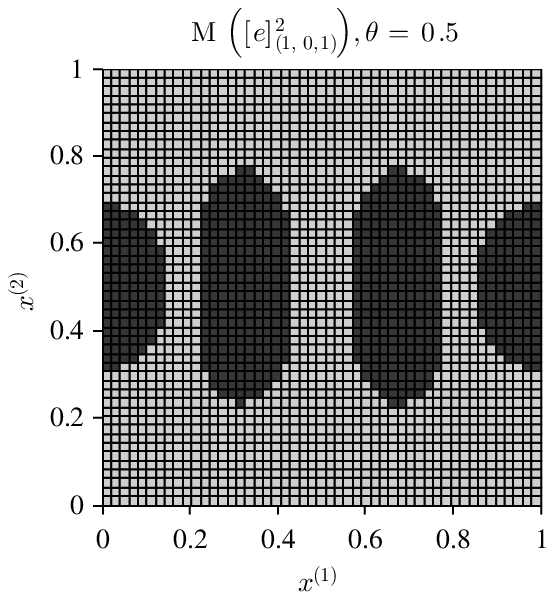}}
	\subfloat[$t^k$ = 4.07]{\includegraphics[trim = 0mm 0mm 0mm 2mm, clip=true, scale=1]{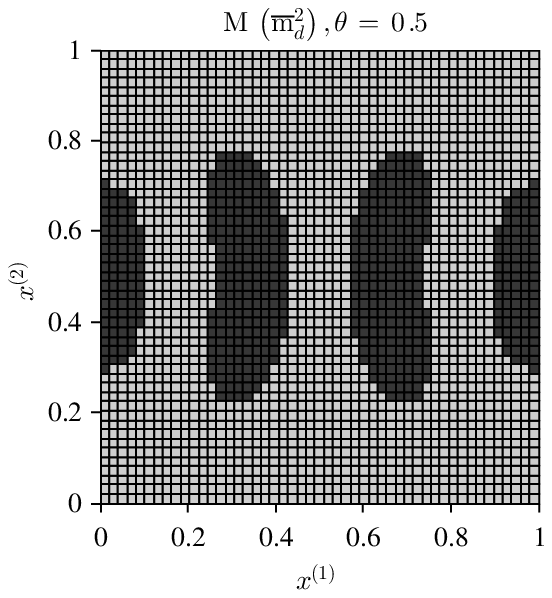}}\\
	\vskip -4pt
	\subfloat[$t^k$ = 8.04]{\includegraphics[trim = 0mm 0mm 0mm 2mm, clip=true, scale=1]{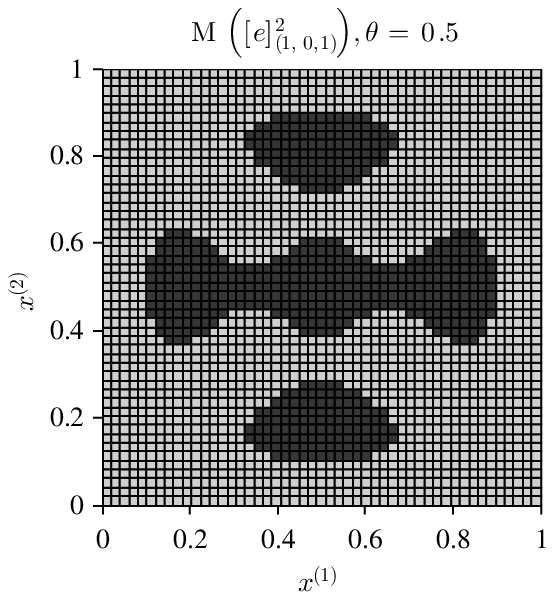}}
	\subfloat[$t^k$ = 8.04]{\includegraphics[trim = 0mm 0mm 0mm 2mm, clip=true, scale=1]{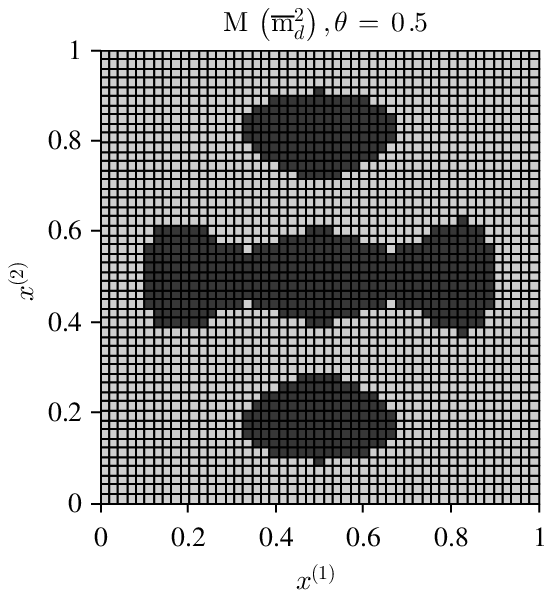}}
	\caption{Example \ref{ex:example-14-2d}.
	'Bulk' marking ($\theta = 0.5$) based on the true error $\error_{(1, \: 0, \: 1)}$
	(left) and indicator $\mdI{} (y = G (\nabla v))$ (right) for $Q^{40}$, $Q^{80}$, and $Q^{80}$.}
	\label{fig:error-maj-bulk-marking-y-ux-ex-14}
\end{figure}
\begin{figure}[!ht]
	\centering
	\subfloat[$t^k$ = 2.06]{\includegraphics[scale=1]{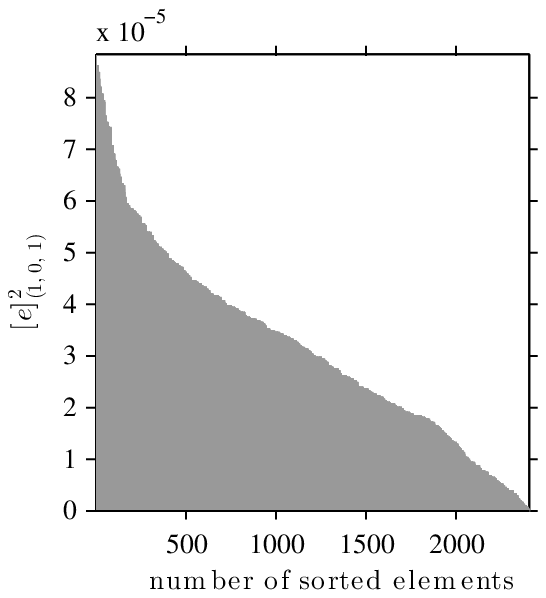}}
	\subfloat[$t^k$ = 2.06]{\includegraphics[scale=1]{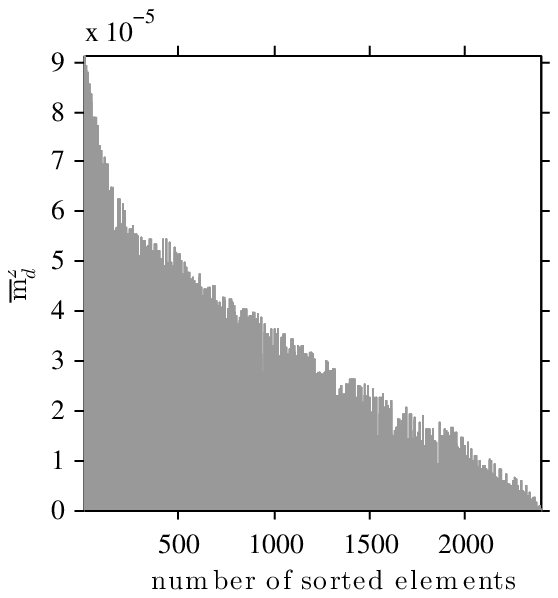}}\\
	\subfloat[$t^k$ = 4.07]{\includegraphics[scale=1]{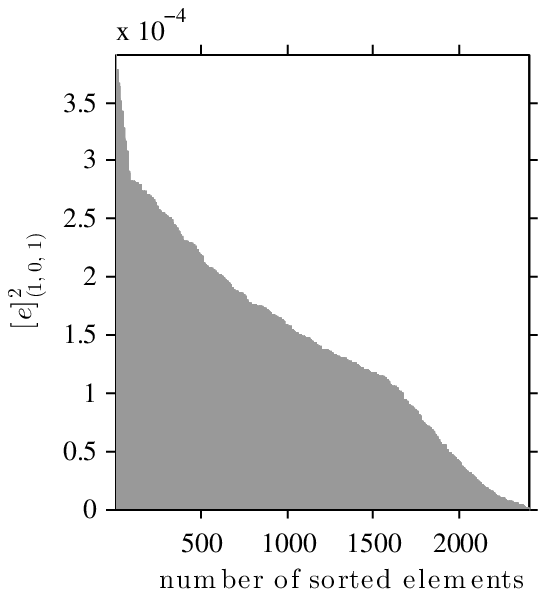}}
	\subfloat[$t^k$ = 4.07]{\includegraphics[scale=1]{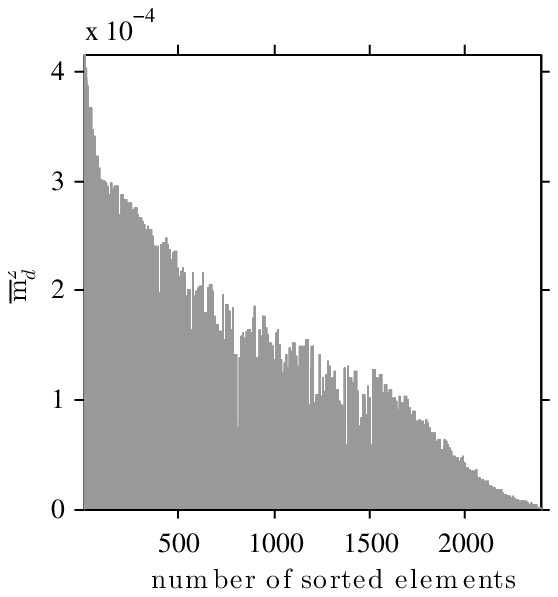}}\\
	\subfloat[$t^k$ = 8.04]{\includegraphics[scale=1]{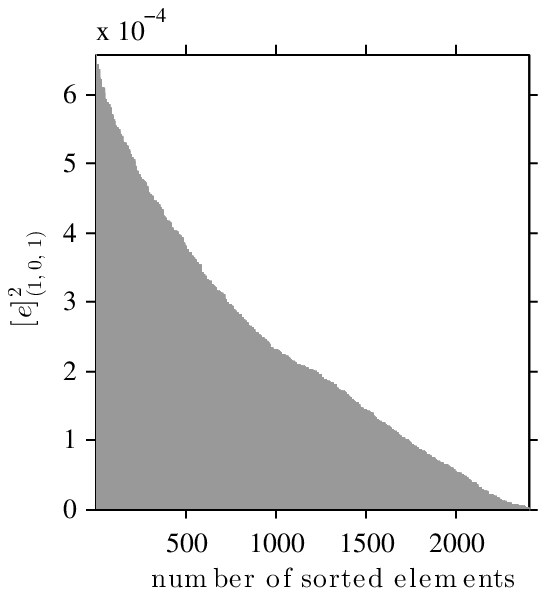}}
	\subfloat[$t^k$ = 8.04]{\includegraphics[scale=1]{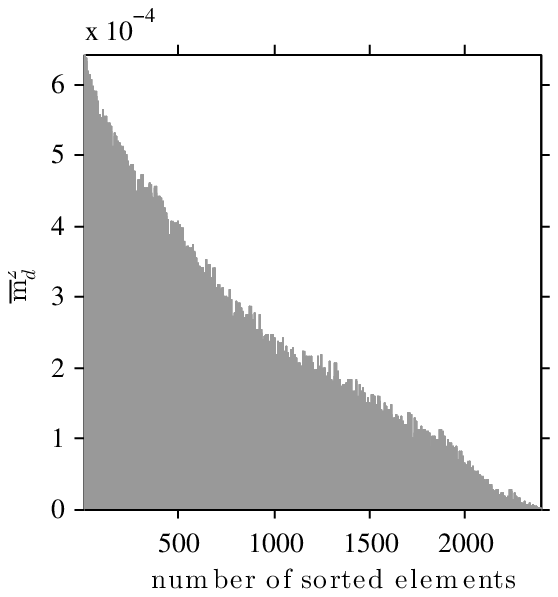}}
	\caption{Example \ref{ex:example-14-2d}.
	Histograms of the ranked element-wise errors and error indicator by $\mdI{}$ for $Q^{40}$, $Q^{80}$, and $Q^{80}$.}
	\label{fig:error-maj-hist-ex-14}
\end{figure}

\end{example}

%\clearpage
%--------------------------------------------------------------------------------------%

\bibliographystyle{abbrv}
\bibliography{main-num}

\begin{thebibliography}{10}

\bibitem{Braess2001}
D.~Braess.
\newblock {\em Finite elements}.
\newblock Cambridge University Press, Cambridge, second edition, 2001.
\newblock Theory, fast solvers, and applications in solid mechanics, Translated
  from the 1992 German edition by Larry L. Schumaker.

\bibitem{Dorfler1996}
W.~D{\"o}rfler.
\newblock A convergent adaptive algorithm for {P}oisson's equation.
\newblock {\em SIAM J. Numer. Anal.}, 33(3):1106--1124, 1996.

\bibitem{Evans2010}
L.~C. Evans.
\newblock {\em Partial differential equations}, volume~19 of {\em Graduate
  Studies in Mathematics}.
\newblock American Mathematical Society, Providence, RI, second edition, 2010.

\bibitem{Johnson2009}
C.~Johnson.
\newblock {\em Numerical solution of partial differential equations by the
  finite element method}.
\newblock Dover Publications Inc., Mineola, NY, 2009.
\newblock Reprint of the 1987 edition.

\bibitem{KleissTomar2012}
S.~Kleiss and S.~Tomar.
\newblock Guaranteed and sharp a posteriori error estimates in isogeometric
  analysis.
\newblock Technical Report~v2, Johann Radon Institute for Computational and
  Applied Mathematics (RICAM), Austrian Academy of Sciences, Altenberger
  Straße 69, A-4040 Linz, Austria, 2013.

\bibitem{Ladyzhenskaya1985}
O.~A. Ladyzhenskaya.
\newblock {\em The boundary value problems of mathematical physics}.
\newblock Springer, New York, 1985.

\bibitem{Ladyzhenskayaetall1967}
O.~A. Ladyzhenskaya, V.~A. Solonnikov, and N.~Uraltseva.
\newblock {\em Linear and quasilinear equations of parabolic type}.
\newblock Nauka, Moscow, 1967.

\bibitem{LeVeque2007}
R.~J. LeVeque.
\newblock {\em Finite difference methods for ordinary and partial differential
  equations}.
\newblock Society for Industrial and Applied Mathematics (SIAM), Philadelphia,
  PA, 2007.
\newblock Steady-state and time-dependent problems.

\bibitem{MaliNeittaanmakiRepin2013}
O.~Mali, P.~Neittaanm{\"a}ki, and S.~Repin.
\newblock {\em Accuracy verification methods. Theory and algorithms (in
  print)}.
\newblock Springer, 2013.

\bibitem{NeittaanmakiRepin2004}
P.~Neittaanm{\"a}ki and S.~Repin.
\newblock {\em Reliable methods for computer simulation}, volume~33 of {\em
  Studies in Mathematics and its Applications}.
\newblock Elsevier Science B.V., Amsterdam, 2004.
\newblock Error control and a posteriori estimates.

\bibitem{NeittaanmakiRepinMaxwell2010}
P.~Neittaanm{\"a}ki and S.~Repin.
\newblock Guaranteed error bounds for conforming approximations of a {M}axwell
  type problem.
\newblock In {\em Applied and numerical partial differential equations},
  volume~15 of {\em Comput. Methods Appl. Sci.}, pages 199--211. Springer, New
  York, 2010.

\bibitem{NeittaanmakiRepin2010}
P.~Neittaanm{\"a}ki and S.~Repin.
\newblock A posteriori error majorants for approximations of the evolutionary
  {S}tokes problem.
\newblock {\em J. Numer. Math.}, 18(2):119--134, 2010.

\bibitem{RepinDeGruyter2008}
S.~Repin.
\newblock {\em A posteriori estimates for partial differential equations},
  volume~4 of {\em Radon Series on Computational and Applied Mathematics}.
\newblock Walter de Gruyter GmbH \& Co. KG, Berlin, 2008.

\bibitem{RepinSauter2006}
S.~Repin and S.~Sauter.
\newblock Functional a posteriori estimates for the reaction-diffusion problem.
\newblock {\em C. R. Acad. Sci. Paris}, 343(1):349--354, 2006.

\bibitem{Repin2002}
S.~I. Repin.
\newblock Estimates of deviations from exact solutions of initial-boundary
  value problem for the heat equation.
\newblock {\em Rend. Mat. Acc. Lincei}, 13(9):121--133, 2002.

\bibitem{RepinSamrowskiSauter2012}
S.~I. Repin, T.~S. Samrowski, and S.~A. Sauter.
\newblock Combined {\it a posteriori} modeling-discretization error estimate
  for elliptic problems with complicated interfaces.
\newblock {\em ESAIM Math. Model. Numer. Anal.}, 46(6):1389--1405, 2012.

\bibitem{RepinTomar2010}
S.~I. Repin and S.~K. Tomar.
\newblock A posteriori error estimates for approximations of evolutionary
  convection-diffusion problems.
\newblock {\em J. Math. Sci. (N. Y.)}, 170(4):554--566, 2010.
\newblock Problems in mathematical analysis. No. 50.

\bibitem{RepinTomar2011}
S.~I. Repin and S.~K. Tomar.
\newblock Guaranteed and robust error bounds for nonconforming approximations
  of elliptic problems.
\newblock {\em IMA J. Numer. Anal.}, 31(2):597--615, 2011.

\bibitem{Thomee2006}
V.~Thom{\'e}e.
\newblock {\em Galerkin finite element methods for parabolic problems},
  volume~25 of {\em Springer Series in Computational Mathematics}.
\newblock Springer-Verlag, Berlin, second edition, 2006.

\bibitem{Valdman2009}
J.~Valdman.
\newblock Minimization of functional majorant in a posteriori error analysis
  based on {$H({\rm div})$} multigrid-preconditioned {CG} method.
\newblock {\em Adv. Numer. Anal.}, pages Art. ID 164519, 15, 2009.

\end{thebibliography}

\end{document}